\documentclass[a4paper,leqno]{amsart}

\usepackage{amsmath}

\usepackage{mathdots} 

\usepackage{amssymb}

\usepackage{amsbsy}

\usepackage{enumerate}

\usepackage[alphabetic]{amsrefs}

\usepackage{mathptmx} 

\usepackage{hyperref}

\usepackage{mathrsfs} 

\usepackage{txfonts}  

\usepackage{stmaryrd} 

\usepackage{fancyhdr}  

\usepackage{comment}

\usepackage{algorithm}

\usepackage{algpseudocode}

\usepackage{tikz}

\usetikzlibrary{arrows,positioning}

\usepackage{makecell}

\usepackage{pgfplots}

\usepackage{graphicx,amssymb}

\usepackage{ifthen}

\usepackage{epsfig}

\newtheorem{theo}{Theorem}

\newtheorem{prop}{Proposition}

\newtheorem{lemm}{Lemma}

\newtheorem{coro}{Corollary}

\newtheorem{defi}{Definition}

\newtheorem{rema}{Remark}

\newtheorem{prob}{Problem}

\newtheorem{exam}{Example}

\renewcommand{\le}{\leqslant}
\renewcommand{\leq}{\leqslant}

\renewcommand{\ge}{\geqslant}
\renewcommand{\geq}{\geqslant}

\newcommand{\la}{\lambda}

\newcommand{\vep}{\varepsilon}

\newcommand{\ds}{\displaystyle}

\newcommand{\mbar}{\ensuremath{\bar{m}}}
\newcommand{\nbar}{\ensuremath{\bar{n}}}

\newcommand{\xbar}{\ensuremath{\bar{x}}}

\newcommand{\bN}{\ensuremath{\mathbb N}}
\newcommand{\bM}{\ensuremath{\mathbb M}}

\newcommand{\bR}{\ensuremath{\mathbb R}}

\newcommand{\bZ}{\ensuremath{\mathbb Z}}

\newcommand{\sI}{\ensuremath{\mathsf I}}

\newcommand{\cA}{\ensuremath{\mathcal A}}
\newcommand{\cB}{\ensuremath{\mathcal B}}
\newcommand{\cC}{\ensuremath{\mathcal C}}

\newcommand{\cQ}{\ensuremath{\mathcal Q}}

\newcommand{\cU}{\ensuremath{\mathcal U}}
\newcommand{\cV}{\ensuremath{\mathcal V}}

\newcommand{\tildef}{\ensuremath{\tilde{f}}}

\newcommand{\dX}{\ensuremath{\mathsf{d_X}}}
\newcommand{\dY}{\ensuremath{\mathsf{d_Y}}}

\newcommand{\sd}{\ensuremath{\mathsf{d}}}

\newcommand{\lip}{\ensuremath{\mathrm{Lip}}}

\newcommand{\diam}{\ensuremath{\mathrm{diam}}}

\newcommand{\met}[1]{\mathsf{#1}}

\newcommand{\banX}{\ban X}
\newcommand{\banY}{\ban Y}
\newcommand{\metX}{\met X}
\newcommand{\metY}{\met Y}
\newcommand{\banXn}{(\ban X,\norm{\cdot})}
\newcommand{\banYn}{(\ban Y,\norm{\cdot})}
\newcommand{\metXd}{(\met X,\dX)}
\newcommand{\metYd}{(\met Y,\dY)}

\newcommand{\eqd}{\stackrel{\mathrm{def}}{=}}

\newcommand{\1}{\mathbf{1}}

\newcommand{\xn}{\ensuremath{\{x_n\}_{n=1}^\infty}}
\newcommand{\yn}{\ensuremath{\{y_n\}_{n=1}^\infty}}

\newcommand{\en}{\ensuremath{(e_n)_{n=1}^\infty}}

\newcommand{\co}{\mathrm{c}_0}

\newcommand{\xbg}{\ensuremath{\left(\sum_{n=1}^\infty\ell_{p_n}\right)_2}}
\newcommand{\xbl}{\ensuremath{\left(\sum_{n=1}^\infty\ell_{\infty}^n\right)_2}}

\newcommand{\ban}[1]{\mathfrak{#1}}

\newcommand{\Nk}{[\mathbb{N}]^{k}}
\newcommand{\Mk}{[\mathbb{M}]^{k}}
\newcommand{\Sk}[2]{{[#1]}^{#2}}

\newcommand{\cube}[1]{\{-1,1\}^{#1}}

\newcommand{\abs}[1]{\lvert #1\rvert}

\newcommand{\norm}[1]{\| #1\|}

\usepackage[alphabetic]{amsrefs}

\begin{document}

\title[Barycentric gluing and geometry of stable metric spaces]{Barycentric gluing and geometry of stable metrics}

\author{F. ~Baudier}
\address{Department of Mathematics, Texas A\&M University, College Station, TX 77843, USA}
\email{florent@math.tamu.edu}
\date{}

\thanks{}
\keywords{}
\subjclass[2010]{46B20, 46B85}

\begin{abstract} We discuss various aspects of a local-to-global embedding technique and the metric geometry of stable metric spaces and of two of its important subclasses: locally finite spaces and proper spaces. We explain how the barycentric gluing technique, which has been mostly applied to bi-Lipschitz embedding problems pertaining to locally finite spaces, can be implemented successfully in a much broader context. For instance, we show that the embeddability of an arbitrary metric space into $\ell_p$ is determined by the embeddability of its balls. We also introduce the notion of upper stability. This new metric invariant lies formally between Maurey-Krivine (metric) notion of stability and Kalton's property $\cQ$. We show that several results of Raynaud and Kalton for stable metrics can be extended to the broader context of upper stable metrics and we point out the relevance of upper stability to a long standing embedding problem raised by Kalton. Applications to compression exponent theory are highlighted and we recall old, and state new, important open problems. This article was written in a style favoring clarity over conciseness in order to make the material appealing, accessible, and reusable to geometers from a variety of backgrounds, and not only to Banach space geometers. 
\end{abstract}

\maketitle

\setcounter{tocdepth}{3}
\tableofcontents
\section{Introduction}
The idea of finding a faithful representation, or embedding, of an intractable space into a tractable host space is an extremely powerful and versatile tool with far reaching applications. The class of spaces of interest and the relevant notions of faithful embeddability are dictated by the problems under scrutiny. This line of thoughts motivated the study of rich and deep embedding theories. Some of the most famous applications of this geometric approach arise in theoretical computer science, geometric group theory, topology, and noncommutative geometry. We refer to \cite{Gromov_ICM84}, \cite{Indyk01}, \cite{Matousek_book02}, \cite{Linial_ICM02}, \cite{Roe_book03}, \cite{Yu_ICM06}, \cite{Naor_ICM10, Naor12, Naor_PNAS13, Naor_ICM18}, \cite{NowakYu_book12}, \cite{Ostrovskii_book13}, \cite{BaudierJohnson16}, \cite{DrutuKapovich_book18} (where comprehensive lists of references can be found) for an extensive account of these rapidly growing research directions. 

The main guiding principle in this article will be to identify the largest possible classes of metric spaces for which every member can be embedded with the highest degree of faithfulness into a Banach space that can be drawn for a class of Banach spaces with the best possible geometric properties. When we think of a Banach space with the best possible geometric properties, Hilbert space is most likely the first space that comes to mind. In most situations (but not all!), Hilbert space exhibits the strongest, or optimal, geometric behavior and we will indeed think of it as a Banach space with highly sought geometric properties. The Lesbegue sequence spaces $\ell_p$ and function spaces $L_p:=L_p[0,1]$ in the super-reflexive range $1<p<\infty$, and more generally super-reflexive Banach spaces\footnote{A Banach space is super-reflexive if all its ultrapowers are reflexive.} are targets of choice due to fundamental applications in non-commutative geometry and topology that will be explained shortly. Despite not possessing the strong geometric features of super-reflexive spaces, $\ell_1$ and $L_1$ nevertheless play a  pivotal role in the design of approximation algorithms. 

We will mainly consider embeddings of metric spaces into Banach spaces and give most definitions for this setting, but obviously these can be generalized to metric space targets. The strongest notion of faithfulness for metric spaces is isometric embeddability. A metric space $(\metX,\dX)$ admits an \emph{isometric embedding} into a Banach space $\banYn$ if there exists $f\colon \metX\to \banY$ such that $\norm{f(x)-f(y)}_\banY=\dX(x,y)$. It follows from a result of Godefroy and Kalton \cite{GodefroyKalton03} that if a Banach space $\banY$ contains an isometric copy of every separable metric space than it must contains a linear isometric copy of every separable Banach space, and thus we cannot expect $\banY$ to have any non-trivial geometric features; in particular $\banY$ cannot be reflexive and even less so super-reflexive. Luckily, for the applications that are of interest to us we only need to consider certain restricted classes of metric spaces and weaker notion of faithfulness. A natural relaxation of isometric embeddability allows for some distortion fo the distances. A \emph{bi-Lipschitz embedding} from $\metXd$ into $\banYn$ is a map $f\colon \metX\to \banY$ satisfying  
\begin{equation}
\dX(x,y)\le \norm{f(x)-f(y)}\leq D\cdot \dX(x,y).
\end{equation}
for all $x,y\in \metX$ and some universal constant $D\in[1,\infty)$. A bi-Lipschitz embedding accounts for the geometry at \emph{all scales} and is the faithfulness notion that is ubiquitous in theoretical computer science, in particular in the design of approximation algorithms. Unlike the situation in the isometric category, it is not true that if a Banach space must contain a bi-Lipschitz copy of every separable metric space then it must contain a linear isomorphic copy of every separable Banach space. Indeed, a landmark result of Aharoni \cite{Aharoni74} states that every separable metric space bi-Lipschitzly embeds into $\co$, but every subspace of $\co$ contains an isomorphic copy of $\co$. Consequently, a Banach space that contains a bi-Lipschitz copy of every separable metric space cannot be super-reflexive, yet even reflexive. Thus, we must once more sacrifice on either the largeness of the class of metric spaces considered or the faithfulness of the embedding. Theoretical computer science is mostly concerned with finite objects, in particular with finite metric spaces, and low-distortion bi-Lipschitz embeddings of finite metric spaces (and subclasses thereof) into finite-dimensional Banach spaces, in particular $\ell_p^k$ for $p\in\{1,2\}$ and a small dimension $k$, has proven to be a crucial tool (c.f.\cite{AIR_ICM18}). We will not expand more on this fascinating area here. We will instead focus our attention to infinite metric spaces and on a relaxation of bi-Lipschitz embeddability, called coarse embeddabilty, which accounts for the large scales only, and 
its small-scale counterpart, namely uniform embeddability. Let $\metXd$ be a metric space, $\banYn$ be a Banach space and $f:\metX\to \banY$ be a map. We define 
\begin{equation*}
\rho_f(t)=\inf\{\norm{f(x)-f(y)}_{\banY} \colon \dX(x,y)\geq t\},
\end{equation*}
and 
\begin{equation*}
\omega_f(t)={\rm sup}\{\norm{f(x)-f(y)}_{\banY}\colon \dX(x,y)\leq t\}.
\end{equation*}
We will refer to the map $f$ as an embedding of $\metX$ into $\banY$ and for every $x,y\in \metX$, it is easy to verify that 
\begin{equation*}
\rho_f(\dX(x,y))\le \norm{f(x)-f(y)}_{\banY}\le\omega_f(\dX(x,y)).
\end{equation*} 

The moduli $\rho_f$ and $\omega_f$ will be called the \emph{compression modulus} and the \emph{expansion modulus} \footnote{the expansion modulus is usually called modulus of uniform continuity when $f$ is a uniformly continuous map.}, respectively, of the embedding. A map $f\colon \metX \to \banY$ is said to be a \emph{uniform embedding} if 
\begin{equation}\label{eq:uniform}
\lim_{t\to 0}\omega_f(t)=0 \text{ and } \rho_f(t)>0 \text{ for all }t>0.
\end{equation} 
This type of embedding is designed to described the microscopic structure of $\metX$ since only the behavior of $f$ with respect to pairs of points whose distance to each other is small is taking into account. It can be seen as a quantitative version of a topological embedding when the spaces carry a metric structure. Uniform embeddability is irrelevant for uniformly discrete spaces and we can always assume that the metric is bounded by replacing $\dX$ with the uniformly equivalent metric $\min\{\dX,1\}$ for instance. Uniform embeddings have been extensively studied in nonlinear Banach space theory for about 50 years, and they recently made a quite unexpected appearance in sketching theory in theoretical computer science (c.f. \cite{AKR18}). 

On the other hand, we say that $\metX$ admits a \emph{coarse embedding}\footnote{For quite some time, in the geometric group theory and noncommutative  geometry communities, a ``coarse embedding'' was simply called a ``uniform embedding'', a shorter version of Gromov's original terminology ``uniform embedding at infinity'' \cite{Gromov93}, but it seems now that the terminology ``coarse embedding'' has been widely adopted.
} into $\banY$ if there is a map $f:\metX\to \banY$ such that 
\begin{equation}\label{eq:coarse}
\lim_{t\to\infty}\rho_f(t)=\infty \text{ and }\omega_f(t)<\infty \text{ for all } t>0.
\end{equation}
It will be interesting to consider embeddings which are \emph{simultaneously} coarse and uniform. A map $f\colon \metX\to \banY$ satisfying \eqref{eq:coarse} and \eqref{eq:uniform} is usually called a \emph{strong embedding}.

The notion of coarse embedding is not relevant for bounded metric spaces, since any map sending all the points of a bounded space onto a single point of the target space is trivially a coarse embedding. In particular, a coarse embedding is not necessarily injective and it is worth noting that we can always assume that $\metX$ is uniformly discrete\footnote{$\metX$ is uniformly discrete if there exists a constant $\delta_s\in(0,\infty)$ such that $\dX(x,y)\ge \delta_s$ for all $x,y\in \metX$. If we want to emphasize on the separation parameter $\delta_s$ we will talk about a $\delta_s$-separated metric space.}  by passing to one of its skeletons. A subset $S$ of a metric space $\metXd$ is called an $(\delta,\Delta)$-\emph{skeleton} if there exist $0<\delta_s\le \delta_m<\infty$ such that $S$ is $\delta_s$-separated and $\sup_{x\in \metX} \dX(x,S)\le \delta_m$. A classical and simple application of Zorn's lemma shows that every non-empty infinite metric space $\metX$ admits a $(\delta,2\delta)$-skeleton for every $\delta\in(0,\diam(\metX))$. The following easy and well-known fact tells us that coarsely embedding a metric space or one of its skeleton is essentially the same, up to some usually inessential loss in the faithfulness of the embedding.
\begin{lemm}\label{lem:skeleton} Let $(\metX,\dX)$ be a metric space and $S$ a $(\delta_s,\delta_m)$-skeleton of $\metX$, then any map $c\colon \metX\to S$ which maps a point in $\metX$ to its closest point in $S$ satisfies for all $x,y\in\metX$
\begin{equation}\label{eq:near-isometry}
\dX(x,y)-2\delta_m\le \dX(c(x),c(y))\le \dX(x,y)+2\delta_m.
\end{equation}
\end{lemm}

The closest point map above is a very faithful coarse embedding\footnote{An embedding satisfying \eqref{eq:near-isometry} is called a \emph{near-isometry} in \cite[Definition 10, page 48]{Rosendal17}}, and it is an example of what is called a \emph{quasi-isometric embedding} in geometric group theory or a \emph{coarse bi-Lipschitz embedding} in nonlinear Banach space theory\footnote{In nonlinear Banach space theory it was customary to say that a space quasi-isometrically embeds if for every $\epsilon>0$ there exists a bi-Lipschitz embedding with distortion at most $1+\epsilon$ but this property is now most commonly referred to as almost isometric embeddability.}.

A map $f\colon \metX\to \banY$ is called a \emph{coarse bi-Lipschitz embedding} or a \emph{quasi-isometric embedding} if it is bi-Lipschitz up to some additive constants, i.e. such that for all $x,y\in \metX$ 
\begin{equation*}
\frac{1}{A}\dX(x,y)-B\le \norm{f(x)-f(y)}_\banY\le A\dX(x,y)+B
\end{equation*}
for some universal constants $A\in[1,\infty)$ and $B\in(0,\infty)$. 

The geometric group theory terminology is well-established and predates the nonlinear Banach space theory terminology. However, coarse bi-Lipschitz embeddings are more tightly connected to bi-Lipschitz embeddings than with isometric embeddings and in the remainder of this article we will favor the terminology ``quasi-isometric'' when groups are involved and the terminology ``coarse bi-Lipschitz'' when treating the case of general metric spaces. Again, the possibility of incorporating a positive constant $B$ (which is not allowed for bi-Lipschitz embedding) allows non injective maps and the small-scale structure can be forgotten in the process.  
In particular, an easy application of Lemma \ref{skeleton} shows that the concepts of coarse bi-Lipschitz embeddability and bi-Lipschitz embeddability for large distances are equivalent, where a map $f\colon \metX\to \banY$ is said to be a \emph{bi-Lipschitz embedding at large distances} if there exist $\tau\in(0,\infty)$ and $D_\tau\in[1,\infty)$ such that $\dX(x,y)\ge\tau$ implies that
\begin{equation*}
\frac{1}{A_\tau}\dX(x,y)\le \|f(x)-f(y)\|_\banY\le D_\tau \dX(x,y).
\end{equation*}

The origins of geometric group theory go back at least as early as Mostow's rigidity theorem and it received a tremendous impetus under Gromov's lead \cite{Gromov93}. For instance, quasi-isometric rigidity of groups has become a prominent branch of geometric group theory \cite{Kapovich14}. The much weaker notion of coarse embeddability is crucial for applications in noncommutative geometry and topology. Gromov suggested in \cite{Gromov95} that a space whose large scale geometry is compatible in a certain sense with the geometry of a super-reflexive Banach space is very likely to satisfy a conjecture of Novikov. Building upon a groundbreaking work of Guoliang Yu \cite{Yu00}, Gromov's intuition was proved to be true by Kasparov and Yu \cite{KasparovYu06}. They showed that if a metric space with bounded geometry\footnote{A metric space has bounded geometry if the number of points in any ball of finite radius is finite and does not depend on the center of the ball, but eventually on the radius.} admits a coarse embedding into a super-reflexive Banach space, then it satisfies the coarse geometric Novikov conjecture. We refer to \cite{FRR95, Valette_book02, XieYu15}, for instance, for more information on the Novikov conjectures and the noncommutative geometry of groups. In the late 90's, the attention was then drawn on the coarse geometry of Banach spaces, which had been little considered at that time, contrary to its uniform counterpart.

In light of Kasparov-Yu result it is natural to ask whether every metric space with bounded geometry admits a coarse embedding into some super-reflexive Banach space. A positive answer to this question would imply the astounding statement that the coarse Novikov conjecture holds for every bounded geometry metric space! This hope will be dashed quickly. As observed by Gromov \cite{Gromov00}, an infinite disjoint union (which has bounded geometry) of a sequence of (finite and regular) expander graphs does not admit a coarse embedding into a Hilbert space. It is much more difficult, but nevertheless true, that there are infinite metric spaces with bounded geometry that do not coarsely embed into \emph{any} super-reflexive Banach space. Two delicate constructions were given by V. Lafforgue \cite{Lafforgue08} (via an algebraic approach), and Mendel and Naor \cite{MN14} (via a graph theoretic approach). However, Brown and Guentner \cite{BrownGuentner05} showed that if we relax the super-reflexivity condition then the situation improves. 

\begin{theo} Every metric space with bounded geometry admits a coarse embedding into the reflexive Banach space $\xbg$, where $\lim_{n\in \bN} p_n=\infty$.
\end{theo}

Unfortunately, the space $\xbg$ is in many aspects far from super-reflexive spaces.  Motivated by the discussion above, constructing embeddings into Banach spaces whose geometric behavior mimics to the greatest extent possible the behavior of super-reflexive spaces has become a fundamental endeavor. Since the appearance of Brown-Guentner embedding result, the theory of metric embeddings of general metric spaces into Banach spaces has sustained a steady growth. An underlying goal of this article if to showcase these developments by emphasizing on the tradeoff between the faithfulness of the embeddings and the strength of the geometric features of the host spaces. 

In Sections \ref{sec:barycentric-lip}-\ref{sec:barycentric-gen} we present a detailed treatment of the barycentric gluing technique introduced in \cite{Baudier07}. In particular we expand its implementation to arbitrary metric spaces and arbitrary embeddings. Rather general new local-to-global embedding results are then obtained, some having already been used in applications (\cite[Theorem 44]{CorderoEskenazis20} for instance). In Sections \ref{sec:Ostrovskii}-\ref{sec:proper-stable}, which are essentially expository, we review top-of-the-art embedding results for the classes of locally finite, proper, and stable metric spaces. Section \ref{sec:upper-stability} focuses on the geometry of stable metric spaces. It is shown that stability, which is an isometric notion, can be relaxed into a bi-Lipschitz variant, called upper stability, that allows us to revisit and extend works of Raynaud and Kalton about stable spaces.

\section{Barycentric gluing: a versatile local-to-global embedding technique}\label{sec:barycentric}

Brown and Guentner  \cite{BrownGuentner05} showed that every bounded geometry metric space admits a coarse embedding into the reflexive Banach space $\xbg$ where $\lim_{n\to \infty} p_n=\infty$. Brown and Guentner embedding technique is a modification of the construction of coarse embeddings using Yu's property $A$ and the compression rate is of the order of $\sqrt{t}$. Using a different embedding technique, Ostrovskii significantly refined Brown and Guentner result when he showed \cite{Ostrovskii06} that every locally finite metric space can be embedded into any Banach space that contains almost isometric copies of $\ell_\infty^n$ for all $n\ge1$ (with a slight abuse of terminology we will simply say ``contains the $\ell_\infty^n$'s'' in the sequel). It worth mentioning that a deep result of Maurey and Pisier \cite{MaureyPisier76} that a Banach space $\banY$ contains the $\ell_\infty^n$'s if and only if $\banY$ does not have finite Rademacher cotype. In particular the host space can be taken to be $\xbg$ but also, and most importantly, $\xbl$. This Banach space is a reflexive Banach space that has the same asymptotic structure as Hilbert space. So from the asymptotic point of view the host Banach space is extremely close to the ``most" super-reflexive Banach space. Moreover, Ostrovskii's embedding compression rate is of the order of\footnote{An inspection of the proof reveals that the exponent could be any number strictly less than~1.} $t^{\log\frac32}$. Whether it is possible to upgrade the previous coarse embeddings to bi-Lipschitz embeddings was unclear at that time and a new idea was needed. 
This idea, that will be referred to as \emph{barycentric gluing}, appeared first in \cite{Baudier07} and was inspired by Ribe's proof \cite{Ribe84} of the uniform equivalence between $\left(\sum_{n=1}^\infty \ell_{q_n}\right)_2$ and $\left(\sum_{n=1}^\infty \ell_{q_n}\right)_2\oplus \ell_1$ where $q_n\to 1$. The original motivation to introduce barycentric gluing was to improve one implication in Bourgain's metric characterization of super-reflexivity in terms of bi-Lipschitz embeddings of finite binary trees \cite{Bourgain86}. Soon after \cite{Baudier07}, it was realized that barycentric gluing was amenable to the study of the bi-Lipschitz geometry of locally finite metric spaces and in \cite{BaudierLancien08} it was shown that Ostrovskii's result holds for the significantly stronger notion of bi-Lipschitz embeddings\footnote{Unfortunately, Ostrovskii's result was not cited in \cite{BaudierLancien08} since the authors were not aware of it at that time. While the result of Brown and Guentner is widely known, Ostrovkii's result seems to have been unintentionaly but unduly overlooked.}. Subsequent implementations of barycentric gluing, which culminated with Ostrovskii's beautiful finite determinacy theorem, relied heavily on the local finiteness of the space to be embedded (\cite{BaudierLancien08}, \cite{Baudier12}, \cite{Ostrovskii12}) or on specific properties of the local embeddings (\cite{BKL10}). As it will be shown in this section, it turns out that barycentric gluing is a much more versatile local-to-global embedding technique than originally perceived. 

\subsection{Embeddability into $L_p$-spaces is locally determined}\label{sec:barycentric-lip}

The following definition is reminiscent of the notion of finite representability in Banach space theory.

\begin{defi}
Let $\lambda\in[1,\infty)$. We will say that a metric space $\metXd$ is \emph{locally $\lambda$-Lipschitz representable} into a metric space $\metYd$ if for every \emph{metric ball} $B$ in $\metX$ there is a map $f_B\colon B\to \met Y$ and a scaling factor $s>0$ such that for all $x,y\in B$ 
\begin{equation}
s\dX(x,y)\le \dY(f_B(x),f_B(y))\le \lambda s\dX(x,y).
\end{equation} 
\end{defi}

In the proof of the next theorem we present the mechanism of the barycentric gluing technique in it most elementary and general implementation. 
 
\begin{theo}\label{thm:barycentric-external}
Let $p\in[1,\infty)$. If $\metXd$ is locally bi-Lipschitzly representable into a Banach space $\banYn$, then $\metX$ bi-Lipschitzly embeds into $\ell_p(\banY)$.
\end{theo}

The proof will actually show that if $\metX$ is locally $\lambda$-Lipschitz representable into $\banY$, then $\metX$ bi-Lipschitzly embeds into $\ell_1(\banY)$ with distortion at most $5\lambda(5\lambda+1)$ and similar estimates can be obtained for other values of $p$.

\begin{proof}

Fix $x_0\in \metX$ and assume, after translating and rescaling if necessary, that for every $k\in\bZ$ there exists $f_k\colon B_\metX(x_0,2^k)\to \banY$ such that $f_k(x_0)=0$ and for all $x,y\in B_\metX(x_0,2^k)$,

\begin{equation}\label{eq:local}
\dX(x,y)\le \|f_k(x)-f_k(y)\|\le \lambda \dX(x,y).
\end{equation}

Note that for all $x\in B_\metX(x_0,2^k)$ one has 

\begin{equation}\label{eq:up-down}
\dX(x_0,x)\stackrel{\eqref{eq:up-down}-\text{low}}{\le} \|f_k(x)\|\stackrel{\eqref{eq:up-down}-\text{up}}{\le} \lambda \dX(x_0,x) \nonumber.
\end{equation} 

Define, $f\colon \metX\to \ell_p(\banY)$ by 

\begin{equation}\label{eq:def}
f(x)=(0,\cdots,0, \mu_x f_k(x),(1-\mu_x)f_{k+1}(x),0,\cdots)
\end{equation}
if $\dX(x_0,x)\in[2^{k-1},2^k)$ where $\mu_x=\frac{2^k-\dX(x_0,x)}{2^{k-1}}$. The map $f$ will actually fall short of providing the desired embedding. Nevertheless, we proceed to estimate the distortion of the embedding $f$ by distinguishing several cases and we will slightly modify $f$ later in order to obtain a genuine bi-Lipschitz embedding. As we shall see it is sufficient to consider the case $p=1$.

Let $x,y\in \metX$ and assume without loss of generality that $\dX(x_0,x)\le \dX(x_0,y)$.

\smallskip

\noindent Observe first that if $k\in\bZ$ and $x\in X$ are such that $2^{k-1}\le \dX(x_0,x)< 2^k$, then
\begin{equation}\label{eq:gen1}
2^k-\dX(x_0,x)\le\mu_x\dX(x_0,x)\le2(2^k-\dX(x_0,x))
\end{equation}
and
\begin{equation}\label{eq:gen2}
\dX(x_0,x)-2^{k-1}\le(1-\mu_x)\dX(x_0,x)\le2(\dX(x_0,x)-2^{k-1})
\end{equation}

\smallskip

\noindent \underline{Case 1.}  $2^{k-1}\le \dX(x_0,x)< 2^k< 2^{r-1}\le \dX(x_0,y)< 2^r$ with $r\ge k+2$. 

\smallskip

Then,
\begin{equation}\label{eq:case1}
\|f(x)-f(y)\|=\mu_x\|f_k(x)\|+(1-\mu_x)\|f_{k+1}(x)\|+\mu_y\|f_r(y)\|+(1-\mu_y)\|f_{r+1}(y)\|,
\end{equation}
and it follows from \eqref{eq:up-down}-up  that

\begin{equation*}
\|f(x)-f(y)\|\le \la (\dX(x_0,x)+\dX(x_0,y)).
\end{equation*}

But in this case, one has $\dX(x_0,x)\le \frac{\dX(x_0,y)}{2}$, therefore

\begin{equation}\label{eq:case1-comp}
\frac{\dX(x_0,y)}{2}\le \dX(x_0,y)-\dX(x_0,x)\le \dX(x,y)\le \dX(x_0,x)+\dX(x_0,y)\le\frac{3\dX(x_0,y)}{2},
\end{equation}

and hence 

\begin{equation}\label{eq:Lip1}
\|f(x)-f(y)\|\le 3\la \dX(x,y).
\end{equation}

The lower bound follows easily from \eqref{eq:case1} and \eqref{eq:up-down}-low

\begin{eqnarray*}
\|f(x)-f(y)\| & \ge & \mu_xd(x_0,x)+(1-\mu_x)d(x_0,x)+\mu_y d(x_0,y)+(1-\mu_y)d(x_0,y)\\
		& \ge & d(x,y).
\end{eqnarray*}

\noindent \underline{Case 2.}  $2^{k-1}\le \dX(x_0,x)< 2^k\le \dX(x_0,y)<2^{k+1}$.

In this case,

\begin{equation}\label{eq:case2}
\|f(x)-f(y)\|=\mu_x\|f_k(x)\|+\|(1-\mu_x)f_{k+1}(x)-\mu_yf_{k+1}(y)\|+(1-\mu_y)\|f_{k+2}(y)\|,
\end{equation}

and it follows from \eqref{eq:up-down}-up and the triangle inequality that

\begin{equation*}\label{Eq:3.11}
\|f(x)-f(y)\|\le\|f_{k+1}(x)-f_{k+1}(y)\|+2\mu_x\la \dX(x_0,x)+2(1-\mu_y)\la \dX(x_0,y).
\end{equation*}

Invoking \eqref{eq:local}, \eqref{eq:gen1} and \eqref{eq:gen2} we obtain

\begin{equation*}\label{Eq:3.12}
\|f(x)-f(y)\|\le \la \dX(x,y)+4\la(2^k-\dX(x_0,x))+4\la(\dX(x_0,y)-2^k).
\end{equation*}

It follows that 

\begin{equation}\label{eq:Lip2}
\|f(x)-f(y)\|\le \la \dX(x,y)+4\la(\dX(x_0,y)-\dX(x_0,x))\le 5\la \dX(x,y).
\end{equation}

For the lower bound, based on \eqref{eq:case2}, \eqref{eq:gen1}, and \ref{eq:gen2} we have  
\begin{equation}\label{eq:aux1}
\|f(x)-f(y)\|\ge \mu_x \dX(x_0,x)\ge 2^k-\dX(x_0,x)
\end{equation}
and 
\begin{equation}\label{eq:aux2}
\|f(x)-f(y)\|\ge (1-\mu_y) \dX(x_0,y)\ge \dX(x_0,y)-2^k
\end{equation}
and summing \eqref{eq:aux1} and \eqref{eq:aux2} gives 
\begin{equation}\label{eq:aux3}
2\|f(x)-f(y)\|\ge \dX(x_0,y)-\dX(x_0,x).
\end{equation}

Note that \eqref{eq:case2} also gives

\begin{align}
\nonumber \|f(x)-f(y)\|&\ge \|f_{k+1}(x)-f_{k+1}(y)\|-\|\mu_xf_{k+1}(x)+(1-\mu_y)f_{k+1}(y)\|\\
\label{eq:aux4}		&\ge \|f_{k+1}(x)-f_{k+1}(y)\|-\mu_x\la \dX(x_0,x)-(1-\mu_y)\la \dX(x_0,y)\\
\label{eq:aux5}		&\ge \dX(x,y)-2\la(2^k-\dX(x_0,x))-2\la(\dX(x_0,y)-2^k)\\
\label{eq:aux6}		&\ge \dX(x,y)-2\la(\dX(x_0,y)-\dX(x_0,x)),
\end{align}
 where we used \eqref{eq:up-down} in \eqref{eq:aux4}, and  \eqref{eq:gen1} and \eqref{eq:gen2} for \eqref{eq:aux5}. Combining \eqref{eq:aux3} and \eqref{eq:aux6} we get 

\begin{equation*}
\|f(x)-f(y)\|\ge \frac{\dX(x,y)}{4\la +1}.
\end{equation*}
\noindent \underline{Case 3.}  $2^{k-1}\le \dX(x_0,x)\leq \dX(x_0,y)< 2^{k}$. 

\smallskip

In this configuration,
\begin{equation}\label{eq:case3}
\|f(x)-f(y)\|=\|\mu_xf_k(x)-\mu_yf_{k}(y)\|+\|(1-\mu_x)f_{k+1}(x)-(1-\mu_y)f_{k+1}(y)\|
\end{equation}
and the triangle inequality provides
\begin{align*}
\|f(x)-f(y)\|\le &\mu_x\|f_k(x)-f_k(y)\|+|\mu_x-\mu_y|\|f_k(y)\|+|\mu_x-\mu_y|\|f_{k+1}(y)\| \\
 &+(1-\mu_x)\|f_{k+1}(x)-f_{k+1}(y)\|.
\end{align*}

Noticing that 
\begin{equation}\label{Eq:3.24}
|\mu_x-\mu_y|=\mu_x-\mu_y=\frac{2^k-\dX(x_0,x)}{2^{k-1}}-\frac{2^k-\dX(x_0,y)}{2^{k-1}}=\frac{\dX(x_0,y)-\dX(x_0,x)}{2^{k-1}}
\end{equation}
it follows from \eqref{Eq:3.24} combined with \eqref{eq:local} and \eqref{eq:up-down}-up that
\begin{align}
\nonumber\|f(x)-f(y)\|&\le \la \dX(x,y)+2\la \dX(x_0,y)\frac{\dX(x_0,y)-\dX(x_0,x)}{2^{k-1}}\\
\label{eq:Lip3}&\le \la \dX(x,y)+4\la (\dX(x_0,y)-\dX(x_0,x))\\
\nonumber	&\le 5\la \dX(x,y).
\end{align}

On the other hand,

\begin{align}
\nonumber\|f(x)-f(y)\|\ge &\mu_x\|f_k(x)-f_k(y)\|+(1-\mu_x)\|f_{k+1}(x)-f_{k+1}(y)\|\\
\nonumber	&-|\mu_x-\mu_y|\|f_k(x)\|-|\mu_x-\mu_y|\|f_{k+1}(y)\|\\
 \label{Eq:3.30}	&\ge \dX(x,y)-2\la(\mu_x-\mu_y) \dX(x_0,y)\\
\nonumber &\ge \dX(x,y)-4\la (\dX(x_0,y)-\dX(x_0,x))
\end{align}

If $\dX(x_0,y)-\dX(x_0,x)\le \frac{\dX(x,y)}{5\la}$ it follows from \eqref{Eq:3.30} that
\begin{equation} \label{Eq:3.32}
\|f(x)-f(y)\|\ge\frac{\dX(x,y)}{5},
\end{equation}

If $\dX(x_0,y)-\dX(x_0,x)\ge \frac{\dX(x,y)}{5\la}$ we cannot conclude. To be able to take care of this inconclusive case we slightly modify the original map $f$ as follows. 

Let $\hat{f}\colon X\to \ell_p(\banY)\oplus_p\bR$ defined by $\hat{f}(x)=(f(x),\dX(x_0,x))$. Since (we are still considering the case $p=1$)
\begin{equation*}
\|\hat{f}(x)-\hat{f}(y)\|=\|f(x)-f(y)\|+|\dX(x_0,x)-\dX(x_0,y)|\le \|f(x)-f(y)\|+\dX(x,y),
\end{equation*}

it follows from \eqref{eq:Lip1}, \eqref{eq:Lip2}, and \eqref{eq:Lip3} that for every $x,y\in X$

\begin{equation*}
\|\hat{f}(x)-\hat{f}(y)\|\le (5\la+1) \dX(x,y).
\end{equation*}

But now if $\dX(x_0,y)-\dX(x_0,x)\ge \frac{\dX(x,y)}{5\la}$ then,
\begin{equation*}
\|\hat{f}(x)-\hat{f}(y)\|\ge |\dX(x_0,x)-\dX(x_0,y)|\ge \frac{\dX(x,y)}{5\la}.
\end{equation*}
Ultimately we proved that there exists a map $\hat{f}\colon X\to \ell_p(\banY)$ such that we have for all $x,y\in X$ (and when $p=1$)
\begin{equation*}
\frac{\dX(x,y)}{5\lambda}\le\|\hat{f}(x)-\hat{f}(y)\|\le (5\lambda+1) \dX(x,y).
\end{equation*} 

To finish the proof in the general case, observe that the Lipschitz constant of $\hat{f}$ is still $5\la+1$ in the case of an arbitrary $p\in[1,\infty)$ since for all $z\in \bR^N$, $\|z\|_p\le \|z\|_1$. As for the lower bound, note that the embedding involves at most 5 coordinates. A careful analysis of the proof shows that the lower bound estimate will only change to $\frac{\dX(x,y)}{5\cdot 2^{1-1/p}}\le\|\hat{f}(x)-\hat{f}(y)\|$ in the range $\la\in[1,2^{1-1/p}]$ (due to the fact that by H\"older's inequality, for all $z\in \bR^N$, $\|z\|_1\le \|z\|_pN^{1-1/p}$) and stay untouched otherwise. 
\end{proof}

It follows from the classical fact that $\ell_p(\ell_p)$ is linearly isometric to $\ell_p$, that the bi-Lipschitz embeddability of a metric space into $\ell_p$ is locally determined.

\begin{coro}\label{coro:1}
Let $p\in[1,\infty)$ and $\metXd$ be a metric space. If $\metX$ is locally Lipschitz-representable into $\ell_p$, then $\metX$ bi-Lipschitzly embeds into $\ell_p$. 
\end{coro}

In Corollary \ref{coro:1}, we could have replaced $\ell_p$ by any Banach space $\banY$ such that $\banY$ is isomorphic to $\ell_p(\banY)$ for some $p$. Theorem \ref{thm:barycentric-external}, and the variety of corollaries that can be deduced from it, is usually sufficient if we are mainly interested in embeddings in the host spaces $\ell_p$ or $L_p[0,1]$. 

In general the host space of the global embedding in Theorem \ref{thm:barycentric-external} will not be isomorphic to the host space of the local embeddings. Nevertheless, in various situations this problem can be overcome. In \cite{Baudier07}, \cite{BaudierLancien08}, \cite{Ostrovskii09}, \cite{BKL10}, and \cite{Baudier12}, barycentric gluing is combined with structural manipulations that can be performed internally to the host space, thus avoiding the use of external $\ell_p$-sums. We now prove a general result where barycentric gluing is implemented internally. Since we will use classical finite-dimensional Schauder decomposition technique, it is crucial that the images of the local embeddings live in finite-dimensional subspaces of the host space. Note that since this will always happen when the domain space is locally finite, it partially explains why barycentric gluing was mainly implemented in the context of locally finite spaces. 

\begin{theo}\label{thm:barycentric-internal}
Let $\metXd$ be a metric space and $\banYn$ a Banach space. Assume that for some $\lambda\ge 1$, 

\begin{enumerate}[(i)]
\item $\metX$ is locally $\lambda$-Lipschitz representable into every finite-codimensional subspace of $\banY$,
\item for every ball in $\metX$, its image under any of the local embeddings from (i) spans a finite-dimensional subspace.
\end{enumerate}
Then, $\metX$ bi-Lipschitzly embeds into $\banY$. 
\end{theo}

\begin{proof}
Fix $x_0\in X$. We will build a finite-dimensional Schauder decomposition inside $Y$ that will be a substitute for the external $\ell_p$-sum in Theorem \ref{thm:barycentric-external}. Let $\gamma >0$ and $B_k\eqd B_X(x_0,2^k)$. Pick a sequence $(\gamma_j)_{j=0}^\infty$ with $\gamma_j>0$ and $\Pi_{j=0}^\infty (1+\gamma_j)\le 1+\gamma$. Choose first a unit vector $v\in \banY$ and let $\bR v$ be the dimension 1 subspace of $\banY$ generated by $v$. Then, using the standard Mazur technique, we can find a finite-codimensional subspace $Z_0$ of $\banY$ so that
\begin{equation*}
\forall y\in \bR v,\ \forall z\in Z_0,\ \|y\|\le (1+\gamma_0)\|y+z\|.
\end{equation*}

By our assumption, and without loss of generality, there exists $f_0\colon B_0\to Z_0$ such that $f_0(x_0)=0$ and for all $x,y\in B_0$,

\begin{equation*}
\dX(x,y)\le \|f_0(x)-f_0(y)\|\le \lambda \dX(x,y).
\end{equation*}

Since by assumption the linear span of the set $f_0(B_0)$, denoted by $H_0$, is a finite-dimensional subspace of $Z_0$, by Mazur's technique we can find one more time a finite-codimensional subspace $Z_1$ of $Y$ so that
\begin{equation*}
\forall y\in \bR v\oplus H_0,\ \forall z\in Z_1,\ \|y\|\le (1+\gamma_1)\|y+z\|.
\end{equation*}
Continuing this process, we can thus construct a sequence $(H_k)_{k=-1}^\infty$ of finite-dimensional subspaces of $Y$ (where $H_{-1}:=\bR v$) and maps $f_k\colon B_k\to H_k$ such that $f_k(x_0)=0$ and for all $x,y\in B_k$,
\begin{equation*}
\dX(x,y)\le \|f_k(x)-f_k(y)\|\le \lambda \dX(x,y).
\end{equation*}
By construction, $(H_j)_{j=-1}^\infty$ is a finite-dimensional Schauder decomposition of its closed linear span $Z$ and there are projections $P_j$ from $Z$ onto $H_{-1}\oplus\dots\oplus H_j$ with kernel $\overline{\rm Span}\,(\bigcup_{i=j+1}^\infty H_i)$ with $\|P_j\|\le 1+\gamma$.
Then we define, $f\colon X\to \oplus_{i=-1}^\infty H_i\subset \banY$ by 
\begin{equation*}
f(x)=d(x_0,x)v\oplus\mu_x f_k(x)\oplus(1-\mu_x)f_{k+1}(x)
\end{equation*}
if $\dX(x_0,x)\in[2^{k-1},2^k)$ where $\mu_x=\frac{2^k-\dX(x_0,x)}{2^{k-1}}$. The proof can be completed by performing computations almost identical to the ones in the proof of Theorem \ref{thm:barycentric-internal} and the details are left to the studious reader. The only difference is that the norm of the projections are not of norm~1 but only bounded above by $1+\gamma$. 
\end{proof}

Even though the assumptions of Theorem \ref{thm:barycentric-internal} might seem quite restrictive, it can be applied successfully in a variety of situations. Observe first that for locally finite metric spaces local Lipschitz representability coincides with the following notion already introduced in \cite[page 1612]{LNP09} (under the terminology \emph{finite representability}).
 
\begin{defi}
Let $\lambda\in[1,\infty)$. We will say that a metric space $\metXd$ is \emph{finitely $\lambda$-Lipschitz representable} into a metric space $\metYd$ if for every \emph{finite subset} $F$ in $\metX$ there is a map $f_F\colon F\to \metY$ and a scaling factor $s>0$ such that for all $x,y\in F$ 
\begin{equation*}
s\cdot \dX(x,y)\le \dY(f_F(x),f_F(y))\le \lambda s\cdot\dX(x,y).
\end{equation*} 
\end{defi}

Therefore, if $\metX$ is locally finite and satisfies assumption $(i)$ in Theorem \ref{thm:barycentric-internal} then it automatically satisfies assumption $(ii)$.
Bourgain \cite{Bourgain86} showed that there exists $\lambda\ge 1$ such that the infinite binary tree is finitely $\lambda$-Lipschitz representable into every non-superreflexive Banach space $X$. Barycentric gluing was originally introduced to show that this implies the bi-Lipschitz embeddability of the infinite binary tree. 

\begin{coro}\label{cor:binary}
The infinite binary tree bi-Lipschitzly embeds into every non-superreflexive Banach space.
\end{coro}

Since super-reflexivity is a finite-codimensional hereditary property, meaning that every co-dimensional subspace of a super-reflexive Banach space is super-reflexive, Corollary \ref{cor:binary}, which originally appeared in \cite{Baudier07}, is an immediate corollary of Theorem \ref{thm:barycentric-internal}.

An immediate consequence of Fr\'echet's embedding \cite{Frechet10} is that every metric space is finitely Lipschitz representable into every Banach which contains the $\ell_\infty^n$'s. It is an easy consequence of Maurey-Pisier theorem that the property of containing the $\ell_\infty^n$'s is a finite co-dimensional hereditary property (see \cite{BaudierLancien15} for a proof) and we can deduce from Theorem \ref{thm:barycentric-internal} the main result of \cite{BaudierLancien08}.


\begin{coro}\label{cor:BL}
Every locally finite metric space bi-Lipschitzly embeds into every Banach space which contains the $\ell_\infty^n$'s. In particular, every localy finite metric space admits a bi-Lipschitzly embedding into the reflexive and asymptotically-$\ell_2$ Banach space $\xbl$.
\end{coro}

Since finite-dimensional subspaces of Hilbert space are simply Euclidean spaces and since every infinite-dimensional Banach space contains the $\ell_2^n$'s by Dvoretzky's theorem \cite{Dvoretzky61}, it follows that Hilbert space is finitely Lipschitz-representable into every infinite-dimensional Banach space.  An immediate consequence of this observation and Theorem \ref{thm:barycentric-internal}, gives the theorem below which first appeared in \cite{Ostrovskii09}.

\begin{coro}\label{cor:Hilbert}
Let $\banY$ be an infinite-dimensional Banach space. Then any locally finite subset of Hilbert space admits a bi-Lipschitz embedding into $\banY$.
\end{coro}

More consequences of Theorem \ref{thm:barycentric-external} to geometric group theory will be discussed in Section \ref{sec:compression}. 

\subsection{Barycentric gluing in the context of general embeddings}\label{sec:barycentric-gen}

In the previous section we focused on bi-Lipschitz embedding. Remarkably, barycentric gluing can be applied to more general notions of embeddings. At this point it is necessary to introduce some terminology. This terminology detour will come in very handy later in order to formulate the results of this section in a condensed fashion. The following definition will be needed in order to quantify the faithfulness of an embedding. 

\begin{defi}[$(\rho,\omega)$-embeddings]
Given non-decreasing functions $\rho,\omega\colon [0,\infty)$ and a map $f\colon\metXd \to \metYd$ such that for all $x,y\in \metX$,
\begin{equation}
\rho(\dX(x,y))\le \dY(f(x),f(y))\le\omega(\dX(x,y))
\end{equation}
we will say that $f$ is a \emph{$[\rho,\omega]$-embedding from $\metX$ into $\met Y$}.
 \end{defi} 
We will mostly focus on coarse and uniform embeddings and compression and expansion rates, and in order to declutter statements from irrelevant (multiplicative) constants we will use a convenient equivalence relation on real functions that capture their large scale and small scale behaviors. For two functions $g,h\colon \bR\to\bR$ we write $g\ll h$ if there exist $c_1,c_2>0$ and $c_3,c_4\in\bR$ such that $g(t)\le c_1h(c_2t+c_3)+c_4$ for every $t\in\bR$. If $g\ll h$ and $h\ll g$ then we write $g\asymp h$. The relation $\asymp$ is easily seen to be an equivalence relation and a function shall be identified with its equivalence class in the sequel. For instance, if we say that $\metX$ admits a $[t/\log (t),t]$-embedding into $\met Y$ this will mean that there exist constants $c_1,c_2,c_3>0$ and $c_4,c_5,c_6,c_7\in \bR$ and a map $f\colon\metXd \to \metYd$ such that for all $x,y\in \metX$,
\begin{equation}
\frac{c_1\dX(x,y)+c_4}{\log(c_2\dX(x,y)+c_5)}+c_6\le \dY(f(x),f(y))\le c_3\dX(x,y)+c_7.
\end{equation}

\begin{rema}\label{rem:skeleton}
Note that Lemma \ref{lem:skeleton} and the convention to identify a function with its equivalence class implies that if a skeleton of $\metXd$ admits a $[\rho,\omega]$-embedding into $\met Y$ then $\metX$ also admits a 
$[\rho,\omega]$-embedding into $\met Y$.
\end{rema}

The notion of Lipschitz representability can naturally be extended in order to accommodate $(\rho,\omega)$-embeddings. 

\begin{defi}[$(\rho,\omega)$-representability]\label{defi:rho-omega}
Let $\rho,\omega\colon [0,\infty)\to [0,\infty)$ be non-decreasing maps. We will say that a metric space $\metXd$ is \emph{locally $(\rho,\omega)$-representable} into a metric space $\metYd$ if there exists for every metric ball $B$ in $\metX$ a map $f_B\colon B\to \metY$ such that for all $x,y\in B$ 
\begin{equation}
\rho(\dX(x,y))\le \dY(f_B(x),f_B(y))\le \omega(\dX(x,y)).
\end{equation} 
\end{defi}

The proof of Theorem \ref{thm:barycentric-external} can be slightly modified to obtain the following more general theorem.

\begin{theo}\label{thm:barycentric-general}
Let $p\in[1,\infty)$. If $\metXd$ is locally $[\rho(t),t]$-representable into a Banach space $\banYn$, then $\metX$ admits a $[\rho(t),t]$-embedding into $\ell_p(\banY)$.
\end{theo}

\begin{proof}
We only emphasize the elements of the proof of Theorem \ref{thm:barycentric-external} that need to be adjusted. Fix $x_0\in \metX$ and without loss of generality assume that for every $k\in\bZ$ there exists $f_k\colon B_\metX(x_0,2^k)\to \banY$ such that $f_k(x_0)=0$ and for all $x,y\in B_\metX(x_0,2^k)$,
\begin{equation*}
\rho(\dX(x,y))-\alpha\le \|f_k(x)-f_k(y)\|\le \lambda\dX(x,y)+\beta,
\end{equation*}
for some $\alpha,\beta\ge 0$ and $\lambda>0$.
Note that this time one has for all $x\in B_\metX(x_0,2^k)$ 
\begin{equation}\label{eq:up-down2}
\rho(\dX(x_0,x))-\alpha\stackrel{\eqref{eq:up-down2}-\textrm{low}}{\le} \|f_k(x)\|\stackrel{\eqref{eq:up-down2}-\textrm{up}}{\le} \lambda\dX(x_0,x)+\beta \nonumber.
\end{equation}
The embedding $\hat{f}\colon X\to \ell_p(\banY)\oplus_p\bR$ is still defined by $\hat{f}(x)=(f(x),\dX(x_0,x))$ where $f\colon \metX\to \ell_p(\banY)$ is again 
\begin{equation*}
f(x)=(0,\cdots,0, \mu_x f_k(x),(1-\mu_x)f_{k+1}(x),0,\cdots)
\end{equation*}
with $\dX(x_0,x)\in[2^{k-1},2^k)$ and $\mu_x=\frac{2^k-\dX(x_0,x)}{2^{k-1}}$. 
In Theorem \ref{thm:barycentric-external} the proof that $\hat{f}$ is $(5\lambda+1)$-Lipschitz only uses \eqref{eq:up-down}-up and substituting \eqref{eq:up-down}-up with 
\eqref{eq:up-down2}-up we can verify that in
\begin{description}
\item[Case 1] $\norm{\hat{f}(x)-\hat{f}(y)}\le (3\lambda+1) \dX(x,y)+2\beta$,
\item[Case 2] $\norm{\hat{f}(x)-\hat{f}(y)}\le (5\lambda+1) \dX(x,y)+5\beta$,
\item[Case 3] $\norm{\hat{f}(x)-\hat{f}(y)}\le (5\lambda+1) \dX(x,y)+3\beta$.
\end{description}

Regarding the lower bound computations (again in the case $p=1$), in \underline{Case 1.} we will now have

\begin{eqnarray*}
\|\hat{f}(x)-\hat{f}(y)\|\ge \|f(x)-f(y)\|\stackrel{\eqref{eq:case1}}{\ge} \mu_y\|f_r(y)\|+(1-\mu_y)\|f_{r+1}(y)\| & \stackrel{\eqref{eq:up-down2}-\text{low}}{\ge} & \rho(\dX(x_0,y))-\alpha\\
 & \stackrel{\eqref{eq:case1-comp}}{\ge} & \rho(\frac{2\dX(x,y)}{3})-\alpha.
\end{eqnarray*}

In  \underline{Case 2.} and \underline{Case 3.}, substituting \eqref{eq:up-down}-low and  \eqref{eq:up-down}-up with  \eqref{eq:up-down2}-low and \eqref{eq:up-down2}-up in the corresponding computations will give 
\begin{equation*}
\|\hat{f}(x)-\hat{f}(y)\|\ge \|f(x)-f(y)\|\ge \rho(\dX(x,y))-\alpha-2(\dX(x_0,y)-\dX(x_0,x)),
\end{equation*}

and 

\begin{equation*}
\|\hat{f}(x)-\hat{f}(y)\|\ge \|f(x)-f(y)\|\ge \rho(\dX(x,y))-\alpha-4(\dX(x_0,y)-\dX(x_0,x)),
\end{equation*}

respectively.

In any case if $\dX(x_0,y)-\dX(x_0,x)\le \frac{\rho(\dX(x,y))}{8}$ it will follow that
\begin{equation*}
\|\hat{f}(x)-\hat{f}(y)\|\ge\frac{\rho(\dX(x,y))}{2}-\alpha,
\end{equation*}

and otherwise

\begin{equation*}
\|\hat{f}(x)-\hat{f}(y)\|\ge |\dX(x_0,x)-\dX(x_0,y)|\ge  \frac{\rho(\dX(x,y))}{8}.
\end{equation*}

Therefore when $p=1$ the map $\hat{f}\colon X\to \ell_p(\banY)$ satisfies for all $x,y\in X$
\begin{equation*}
\frac{\rho(2\dX(x,y)/3)}{8}-\alpha\le\|\hat{f}(x)-\hat{f}(y)\|\le (5\lambda+1) \dX(x,y)+5\beta.
\end{equation*}

The observations regarding the validity of the proof for an arbitrary $p\in[1,\infty]$  are still valid and the conclusion follows.

\end{proof}

Notice that in the proof of Theorem \ref{thm:barycentric-general} the only property of $\rho$ that we use is its monotonicity and thus Theorem \ref{thm:barycentric-general} can be applied in many situations. 
As a typical example, we wish to spell out a particularly interesting application of Theorem \ref{thm:barycentric-general} to coarse geometry.

Recall first that a family of metric spaces $(\metX_i)_{i\in I}$, such that $\sup_{i\in I}\diam(\metX_i)=\infty$, is said to be \emph{equi-coarsely embeddable} into $\metYd$ if there exist $\rho, \omega\colon [0,\infty)\to[0,\infty)$ non-decreasing with $\lim_{t\to\infty}\rho(t)=\infty$ and for every $i\in I$ a map $f_i\colon \metX_i\to \met Y$ such that for all $x,y\in \metX_i$ $$\rho(\sd_{\metX_i}(x,y))\le \dY(f_i(x),f_i(y))\le\omega(\dX_i(x,y)).$$

Note that by Remark \ref{rem:skeleton}, the proof of Theorem \ref{thm:barycentric-general} will work just fine if if merely assume that for some $x_0\in \metX$ the sequence of dyadic balls $\{B_\metX(x_0,2^n)\}_{n\ge 1}$ admits an equi-coarse embedding into a Banach space $\banY$ with compression function $\rho$ and expansion function $\omega(t)\asymp t$. 
 
 \begin{coro}
 Let $\metXd$ be a metric space. Assume that for some $x_0\in \metX$ the sequence of dyadic balls $\{B_\metX(x_0,2^n)\}_{n\ge 1}$ admits an equi-coarse embedding into a Banach space $\banY$ with compression function $\rho$ and expansion function $\omega(t)\asymp t$, then for any $p\in[1,\infty]$ $\metX$ admits a coarse embedding into $\ell_p(\banY)$ with compression function $\tilde{\rho}\asymp \rho$ and expansion function $\tilde{\omega}(t)\asymp t$.
 \end{coro}

\begin{rema}
For metrically convex spaces one can always assume that $\omega(t)\asymp t$, and thus for coarse embeddings of metrically convex spaces we only need to focus on the compression function $\rho$.
\end{rema}

The proof of Theorem \ref{thm:barycentric-general} can also be adjusted to obtain.

\begin{theo}
Let $\metXd$ be a metric space and $\banYn$ a Banach space. Assume that 
\begin{enumerate}[(i)]
\item $\metX$ is locally $[\rho(t),t]$-representable into every finite-codimensional subspace of $\banY$,
\item for every ball in $\metX$, its image under any of the local embeddings from (i) spans a finite-dimensional subspace.
\end{enumerate}
Then, $\metX$ admits a $[\rho(t),t]$-embedding into $\banY$. 
\end{theo}
 
\subsection{Ostrovskii's finite determinacy theorem}\label{sec:Ostrovskii}
In Section \ref{sec:barycentric-lip} we saw that an arbitrary metric space admits a bi-Lipschitz embedding into $\ell_p$ whenever it is locally Lipschitz representable in $\ell_p$. For arbitrary Banach space targets, the same conclusion holds under the stronger representability requirements of Theorem \ref{thm:barycentric-general}. If we restrict ourselves to the much smaller class of locally finite metric spaces, the representability requirements of Theorem \ref{thm:barycentric-general} reduce to finite Lipschitz representability into \emph{every} finite-codimensional subspace. In 2012, Ostrovskii showed that this stronger representability requirement is superfluous. 

\begin{theo}\label{thm:OFD}
If a locally finite metric space $\metXd$ is finitely Lipschitz representable into an arbitrary Banach space $\banYn$ then $\metX$ admits a bi-Lipschitz embedding into $\banY$.
\end{theo}

Theorem \ref{thm:OFD} will be referred to as \emph{Ostrovskii's finite determinacy theorem}, and we will say that the bi-Lipschitz embeddability of locally finite metric spaces into Banach spaces is \emph{finitely determined}.

Recall that a Banach space $\banX$ is said to be \emph{finitely $\lambda$-representable} in another Banach space $\banY$ if for every finite-dimensional subspace $F$ of $\banX$ there exists a finite-dimensional subspace $G$ of $\banY$ and a bounded linear map $T\colon F\to G$ such that $\dim(F)=\dim(G)$ and $\|T\|\cdot\|T^{-1}\|\le \lambda$. We then say that $X$ is \emph{finitely representable} in $Y$ if it is $\lambda$-crudely finitely representable for some $\lambda\in[1,\infty)$. The heart of the proof of Ostrovskii's finite determinacy theorem is the following result (which is implicit in \cite{Ostrovskii12}). 

\begin{theo}\label{thm:cfrLipschitz}
If a Banach space $\banX$ is finitely representable into a Banach space $\banY$ then every locally finite subset of $\banX$ admits a bi-Lipschitz embedding into $\banY$. 
\end{theo}

The assumption of Theorem \ref{thm:cfrLipschitz} implies that every locally finite subset $\met{M}$ of $\banX$ is finitely $\lambda$-Lipschitz representable into $\banY$, for some $\lambda\ge 1$, and thus an external application of barycentric gluing (Theorem \ref{thm:barycentric-external} obviously applies here) gives that $\met{M}$ bi-Lipschitzly embeds into $\ell_p(\banY)$. The key new idea of Ostrovskii was to show that, for locally finite metric spaces, it is possible to implement the barycentric gluing internally via a delicate compactness argument. The proof of Theorem \ref{thm:cfrLipschitz} has many technical intricacies and we will only sketch the main lines of Ostrovskii's argument.

\begin{proof}[Sketch of proof of Theorem \ref{thm:cfrLipschitz}]
Let $\met{M}$ be a locally finite subset of $\banX$. By translating and rescaling if necessary we can assume without loss of generality that $0\in \met{M}$ and that $\norm{x}_\banX\ge 1$ for all $x\in\met{M}\setminus\{0\}$. It follows easily form the assumption of Theorem \ref{thm:cfrLipschitz} that for some $\lambda \ge 1$ such that for all $k\ge 1$ there exists a map $f_k \colon B_k\eqd\{x\in \met{M}; \|x\|_\banX\le 2^{k}\} \to \banY$ such that $f_k(0)=0$ and for all $x,y\in B_k$
\begin{equation}\label{fn}
\|x-y\|_\banX\le\|f_k(x)-f_k(y)\|_\banY\le\lambda \|x-y\|_\banX.
\end{equation}
We can then naturally define $f(x)=\mu_x f_{k}(x)+(1-\mu_x)f_{k+1}(x)$, if $2^{k-1}\le\| x\|_\banX\le 2^{k}$ where $\mu_x=\frac{2^{k}-\| x\|_\banX}{2^{k-1}}$. Essentially the same computations as in the proof of Theorem \ref{thm:barycentric-external} will show that $f$ is a Lipschitz map.
To handle the case where $x$ and $y$ lie in the same dyadic annulus in the proof of Theorem \ref{thm:barycentric-external}, we had to slightly modify the embedding. We also need to do this here and we could modify $f$ similarly and define $\tilde{f}(x)=(f(x),\|x\|)\in Y\oplus \bR$. This would be fine if $Y\oplus \bR$ were isomorphic to $Y$. Unfortunately, the existence of Banach spaces that are not isomorphic to any of their hyperplans is an annoyance that can be dealt with the following lemma (see \cite{Ostrovskii12} for a proof).

\begin{lemm}\label{lem:modif}
There exist a $1$-Lipschitz map $\tau\colon[0,\infty)\to \banY$ and $c\in(0,\infty)$ such that if $\tildef\colon \met{M}\to \banY$ is given by $\tildef(x)=\tau(\|x\|)+f(x)$ then $\tildef$ is Lipschitz and for every $x,y\in \met{M}$, we have $\|\tildef(x)-\tildef(y)\|\ge c(\|y\|-\|x\|)$.
\end{lemm}

It is easy to verify that if there is at least a dyadic annulus separating $x$ and $y$, i.e. $2^{k-1}\leq\|x\|<2^{k}<2^r\leq \|y\|<2^{r+1}$, for some $k,r\in\bN$, with $r\ge k+2$, then $\|y\|-\|x\|\ge\frac{\|x-y\|}{4}$ and hence the modification $\tildef$ is co-Lipschitz. It is also easy to verify, using only the triangle inequality, that if $x$ and $y$ are in two consecutive annuli then $\|f(x)-f(y)\|\ge\|x-y\|-8\lambda(\|y\|-\|x\|)$. If $\|y\|-\|x\|$ is small compared to $\|x-y\|$ then $f$ (and thus $\tilde{f}$) is already co-Lipschitz, otherwise the modification $\tilde{f}$ will be co-Lipschitz. 

The most delicate configuration is, as in Theorem \ref{thm:barycentric-external}, when the two points are in the same dyadic annulus. In this case, we do not have enough information on the $f_k$'s to estimate the co-Lipschitz constant of $\tilde{f}$. Ostrovskii's solution was to use auxiliary accumulation points of the embedding to perform the lower bound computations. More precisely, for every $x\in \met{M}$ the sequence $\{f_k(x)\}_{k\ge 1}$ is bounded (note that $f_k(x)$ is only defined for all $k$ larger that a certain value, but this is inessential). We will often conveniently identify a vector in $\banY$ with its canonical image in the bidual $\banY^{**}$. Note that the bounded subsets of $\banY^{**}$ are weak*-compact but not necessarily sequentially  weak*-compact since $\banY^*$ might not be separable. However since $\met{M}$ is locally finite we can argue that the union of the images of the local embeddings can be identify with a subset of a Banach space $Z^{**}\subset \banY^{**}$ such that $Z^*$ is separable and without loss of generality we will assume that $\banY$ and $\banY^*$ are separable. Since $\met{M}$ is locally finite, it follows from a classical diagonal argument that we can pick a subsequence $\{f_{k_n}\}_{n\ge 1}$ such that $f_{k_n}(x)\stackrel{w*}{\to} f_\omega(x)\in \banY^{**}$ for all $x\in \met{M}$. It is important to note that since the sequence of sets $\{B_k\}_{k=1}^\infty$ is non-decreasing, for any subsequence $\{f_{k_n}\}_{n=1}^{\infty}$ of $\{f_{k}\}_{k=1}^{\infty}$,  $f_{k_n}$ will map $B_{k}$ into $\banY$ for all $k\le n$, and will satisfy $(\ref{fn})$. Therefore, the arguments given above will still be valid for any subsequence of $\{f_{k}\}_{k=1}^{\infty}$. 

It remains to carefully select a subsequence so that we can estimate the co-Lipschitz constant in the case where the pair of points lies in the same dyadic annulus. We briefly describe how this selection process is done in \cite{Ostrovskii12}. By Hahn-Banach theorem, for every vector $y\in \banY$ there is a linear functional $g\in \banY^*$ of norm~1 such that $g(y)=\|y\|$. If a sequence $\yn$ of vectors in $\banY$ weak*-converges to $z\in \banY^{**}$, then for any linear functional $g\in \banY^*$, we have $\lim_n (y_n-z)(g)=0$. In the sequel we will make a convenient abuse of notation and write $g(z)$ for $z(g)$ when $g\in \banY^*$ and $z\in \banY^{**}$. For some technical reasons we want the maps $f_n$ to satisfy some additional properties. Extracting as in \cite{Ostrovskii12} we can find a collection of vectors $(f_\omega(x))_{x\in \banX}$ and a collection of linear forms such that : 
\begin{itemize}
\item If for some $x,y\in B_k$ and some $m\ge k$ the vector $f_m(x)-f_m(y)-(f_\omega(x)-f_\omega(y))$ is nonzero, there exists $g_{m,x,y}\in \banY^*$ of norm $1$ such that
\begin{equation}\label{ineq1}
g_{m,x,y}(f_m(x)-f_m(y)-(f_\omega(x)-f_\omega(y))\ge\frac{99}{100}\|f_m(x)-f_m(y)-(f_\omega(x)-f_\omega(y))\|,
\end{equation}
and 
\begin{equation}\label{ineq2}
\textrm{for all } p>m,\, |g_{m,x,y}(f_p(x)-f_p(y)-(f_\omega(x)-f_\omega(y))|\le\frac{1}{1000}\|x-y\|.
\end{equation}
\item If $f_\omega(x)\neq f_\omega(y)$ for $x,y\in B_k$, there exists $h_{x,y}\in \banY^*$ of norm $1$ such that
\begin{equation}\label{ineq3}
h_{x,y}(f_\omega(x)-f_\omega(y))\ge \frac{99}{100}\|f_\omega(x)-f_\omega(y)\|
\end{equation}
and
\begin{equation}\label{ineq4}
\textrm{for all } m\ge k,\, |h_{x,y}(f_m(x)-f_m(y)-(f_\omega(x)-f_\omega(y))|\le\frac{1}{100}\|f_\omega(x)-f_\omega(y)\|.
\end{equation}
\end{itemize}

Expressions of the form $\|\tildef(x)-\tildef(y)\|$ can now be estimated from below with the help of the linear forms of norm $1$ satisfying inequalities (\ref{ineq1}) and (\ref{ineq2}), or (\ref{ineq3}) and (\ref{ineq4}).

\smallskip

\noindent \underline{Case 1.} $2^{k-1}\leq \|x\|\leq \|y\|<2^{k}$, for some $k\in\bN$.

\smallskip

It is clear from the property of $\tildef$ that there is nothing to prove if $\|y\|-\|x\|\ge\frac{\|x-y\|}{1000\lambda}$. From now on we assume that $\|y\|-\|x\|\le\frac{\|x-y\|}{1000\lambda}$.

\noindent \underline{Subcase 1.a.} $\|f_{\omega}(x)-f_{\omega}(y)\|\ge\frac{1}{100}\|x-y\|$ 
 \begin{align*}
 f(x)-f(y)=&f_{\omega}(x)-f_{\omega}(y)+\mu_x(f_{k}(x)-f_{k}(y)-(f_{\omega}(x)-f_{\omega}(y)))\\
               &+(1-\mu_x)(f_{k+1}(x)-f_{k+1}(y)-(f_{\omega}(x)-f_{\omega}(y)))\\
               &+(\mu_x-\mu_y)f_{k}(y)+(\mu_y-\mu_x)f_{k+1}(y)
 \end{align*}  
 and hence
  \begin{eqnarray*}
 \|f(x)-f(y)\| & \ge &h_{x,y}(f(x)-f(y)) \\
	  & \ge & h_{x,y}(f_{\omega}(x)-f_{\omega}(y))+\mu_xh_{x,y}\left(f_{k}(x)-f_{k}(y)-(f_{\omega}(x)-f_{\omega}(y))\right)\\
	& & +(1-\mu_x)h_{x,y}(f_{k+1}(x)-f_{k+1}(y)-(f_{\omega}(x)-f_{\omega}(y)))\\
	& & +(\mu_x-\mu_y)h_{x,y}(f_{k}(y))+(\mu_y-\mu_x)h_{x,y}(f_{k+1}(y))\\
	  & \ge & \frac{99}{100}\|f_{\omega}(x)-f_{\omega}(y)\|-\mu_x\frac{1}{100}\|f_{\omega}(x)-f_{\omega}(y)\|\\
	  & & -(1-\mu_x)\frac{1}{100}\|f_{\omega}(x)-f_{\omega}(y)\|-2|\mu_x-\mu_y|\lambda\|y\|\\
           & \ge & \frac{98}{100}\|f_{\omega}(x)-f_{\omega}(y)\|-4\lambda(\|y\|-\|x\|)\\
           & \ge & \frac{98}{10000}\|x-y\|-\frac{4}{1000}\|x-y\|\\
           &  \ge & \frac{48}{10000}\|x-y\|.
 \end{eqnarray*}  

\noindent \underline{Subcase 1.b.} $\|f_{\omega}(x)-f_{\omega}(y)\|<\frac{1}{100}\|x-y\|$

\noindent \underline{Subcase 1.b.i.} $\mu_x\|f_{k}(x)-f_{k}(y)-(f_{\omega}(x)-f_{\omega}(y))\|\ge\frac{1}{10}\|x-y\|$
                   
  \begin{eqnarray*}
 \|f(x)-f(y)\| & \ge & g_{k,x,y}(f(x)-f(y))\\
                 & \ge & \mu_xg_{k,x,y}(f_{k}(x)-f_{k}(y)-(f_{\omega}(x)-f_{\omega}(y)))+g_{k,x,y}(f_{\omega}(x)-f_{\omega}(y))\\
                 & & +(1-\mu_x)g_{k,x,y}(f_{k+1}(x)-f_{k+1}(y)-(f_{\omega}(x)-f_{\omega}(y)))\\
               & & +(\mu_x-\mu_y)g_{k,x,y}(f_{k}(y))+(\mu_y-\mu_x)g_{k,x,y}(f_{k+1}(y))\\
               & \ge & \frac{99}{1000}\|x-y\|-\frac{1}{100}\|x-y\|-\frac{1}{1000}\|x-y\|-\frac{4}{1000}\|x-y\|\\
               & \ge & \frac{84}{1000}\|x-y\|.
 \end{eqnarray*} 
  
\noindent \underline{Subcase 1.b.ii.} $\mu_x\|f_{k}(x)-f_{k}(y)-(f_{\omega}(x)-f_{\omega}(y))\|<\frac{1}{10}\|x-y\|$ 

Remark that for any $n$, in particular for $n=k$ or $k+1$, one has
\begin{align*}
\|f_{n}(x)-f_{n}(y)-(f_{\omega}(x)-f_{\omega}(y))\|&\ge \|f_{n}(x)-f_{n}(y)\|-\|f_{\omega}(x)-f_{\omega}(y)\|\\
                                                &\ge \|x-y\|-\frac{\|x-y\|}{100}\\
                                                &\ge\frac{99}{100}\|x-y\|,
\end{align*}

and hence $\mu_x<\frac{10}{99}$ and $1-\mu_x>\frac{89}{99}$.  Consequently,

 \begin{eqnarray*}
 \|f(x)-f(y)\| & \ge & (1-\mu_x)\|f_{k+1}(x)-f_{k+1}(y)-(f_{\omega}(x)-f_{\omega}(y))\|\\
               & &-\|f_{\omega}(x)-f_{\omega}(y)\|-\mu_x\|f_{k}(x)-f_{k}(y)-(f_{\omega}(x)-f_{\omega}(y))\|\\
               & &-\|(\mu_x-\mu_y)f_{k}(y)+(\mu_y-\mu_x)f_{k+1}(y)\|\\
               & \ge & \frac{89}{99}\frac{99}{100}\|x-y\|-\frac{1}{100}\|x-y\|-\frac{1}{10}\|x-y\|-\frac{4}{1000}\|x-y\|\\
               &\ge &\frac{776}{1000}\|x-y\|
 \end{eqnarray*}  

\end{proof}

It is well known that any ultrapower of a Banach space $\banX$ is finitely representable in the space $\banX$ itself.

\begin{coro}\label{coro:ultra} Let $\banX$ be a Banach space, and $\cU$ be any non-principal ultrafilter, then every locally finite subset of $\banX^{\cU}$ admits a bi-Lipschitz embedding into $\banX$.
\end{coro}

Theorem \ref{thm:OFD} now follows easily from Corollary \ref{coro:ultra}.

\begin{proof}[Proof of Theorem \ref{thm:OFD}]
Let $\metX$ be locally finite, and $\banY$ be a Banach space. Fix $x_0\in \metX$ and assume that $\metX$ is finitely $\lambda$-Lipschitz representable into $Y$. Then for some $x_0\in \metX$ and for every $k\in\bN$ there exists $f_k\colon B(x_0,2^k)\to \banY$ such that $f_k(x_0)=0$ and for every $x,y\in B(x_0,2^k)$, 
\begin{equation*}
\dX(x,y)\le\|f_k(x)-f_k(y)\|_\banY\le \lambda \dX(x,y).
\end{equation*}
For all $x\in \metX$, the sequence $(f_{k}(x))_k$ is bounded and we define 
$$\begin{array}{rcl}
    f\colon \metX & \rightarrow & \banY^\cU \\
          x & \mapsto & (0,\dots,0,f_{k}(x), f_{k+1}(x),\dots),\textrm{ if } 2^{k-1}\le \dX(x,x_0)\le 2^{k}\textrm{ for some }k\in\bN.\\
    \end{array}$$
 By definition of the norm of the ultrapower it is clear that for every $x,y\in \metX$, 
 \begin{equation*}
 \dX(x,y)\le\|f(x)-f(y)\|_{\banY^\cU}\le \lambda \dX(x,y).
 \end{equation*}
Since, local finiteness is preserved under bi-Lipschitz embeddings,  $f(\metX)$ is a locally finite subset of $Y^{\cU}$. The composition of a bi-Lipschitz embedding and another bi-Lipschitz embedding being a bi-Lipschitz embedding, the conclusion follows from Corollary \ref{coro:ultra}.
\end{proof}

While the statement of Theorem \ref{thm:OFD} is conceptually very appealing, we feel that Theorem \ref{thm:cfrLipschitz} and Corollary \ref{coro:ultra} should be more useful in practice. For example, the next corollary, which seems new, can be deduced from Corollary \ref{coro:ultra}. The equivalence in Corollary \ref{cor:equivqilip} is an analogue in the locally finite context, of a well-known phenomena for families of unbounded finite metric spaces.

\begin{coro}\label{cor:equivqilip}
Let $\metX$ be a locally finite metric space, and $\banY$ be a Banach space then, $\metX$ admits a coarse bi-Lipschitz embedding into $\banY$ if and only if $\metX$ admits a bi-Lipschitz embedding into $\banY$. 
\end{coro}

\begin{proof}
One implication is trivial. Assume without loss of generality assume that $\metX$ is a locally finite subset of the Banach space $\co$, and that there exists $f\colon\metX\subset \co\to \banY$ such that $$\frac{1}{A}\|x-y\|_{\infty}-B\le\|f(x)-f(y)\|_\banY\le A\|x-y\|_{\infty}+B.$$ Define
$$\begin{array}{rcl}

    f^{\cU}\colon \metX & \rightarrow & \banY^\cU \\
            &   &  \\
          x & \mapsto & (\frac{f(nx)}{n})_{n\ge 1},\\
    \end{array}$$
then $\frac{1}{A}\|x-y\|_{\infty}\le\|f^{\cU}(x)-f^{\cU}(y)\|_{\banY^{\cU}}\le A\|x-y\|_{\infty}$, for all $x,y\in\metX$. Since obviously $f^{\cU}(\metX)$ is a locally finite subset of $\banY^{\cU}$, by Corollary \ref{coro:ultra} there exists a bi-Lipschitz map $g\colon f^{\cU}(\metX)\subset \banY^{\cU}\to \banY$ and $g\circ f^{\cU}$ is a bi-Lipschitz embedding of $\metX$ into $\banY$.
\end{proof}

It is clear from the proof above that the embeddability of locally finite metric spaces is finitely determined for other notions of embeddability, as long as: 
\begin{enumerate}
\item the composition with a bi-Lipschitz embedding does not change the nature of the embedding,
\item the induced notion of representability requires a uniform control on the expansion and compression moduli,
\item the embedding preserves local finiteness.
\end{enumerate} 

\begin{coro}\label{coro:coarsefindet}
The coarse embeddability of a locally finite metric space into an arbitrary Banach space is finitely determined.
\end{coro}

\section{Simultaneously coarse and uniform embeddings of proper and stable metrics}\label{sec:proper-stable}
In the previous section we saw that every locally finite metric space (and thus every metric space with bounded geometry) admits a bi-Lipschitz embedding into every Banach space that contains the $\ell_\infty^n$'s. It is natural to investigate if similar general embedding results hold for larger classes of metric spaces. In this section, which is of expository nature, we discuss satisfactory answers that have been given for proper metrics in \cite{Baudier12,BaudierLancien15} and for stable metrics in \cite{Kalton07}. 

\subsection{Almost bi-Lipschitz embeddability of proper metric spaces}\label{sec:proper}
A metric space is \emph{proper} if all its closed balls are compact. It is plain that locally finite metrics are proper metrics and that the class of proper metrics strictly contains the class of locally finite metrics. Observing that any skeleton of a proper metric space is locally finite and recalling that every skeleton is coarse-Lipschitz equivalent to the original space, it follows that every embedding result in the previous section pertaining to locally finite metrics extends to an embedding result for proper metrics.  We just single out one such example and leave other similar easy derivations to the reader.  

\begin{coro}\label{cor:proper-CL}
Every proper metric space coarse bi-Lipschitzly embeds into every Banach space with trivial cotype. In particular, every proper metric space admits a coarse bi-Lipschitzly embedding into the reflexive and asymptotically-$\ell_2$ Banach space $\xbl$.
\end{coro}

An interesting problem for proper metric spaces is to preserve the small-scale structure (via uniform embeddings) or, more ambitiously, to preserve simultaneously the small-scale and the large-scale structure (via embeddings that are simultaneously uniform and coarse). Consider first the case of compact metric spaces. A compact metric space being bounded it does not make much sense to discuss its large-scale geometry. However, as far as its small-scale geometry is concerned, it was known since \cite[Proposition 7.18]{BenyaminiLindenstrauss00} that for every compact metric space $(\met{X},\dX)$ there exist a constant $C>0$ and a map $f\colon \metX \to\xbl$ such that for all $x,y\in \metX$, 
\begin{equation}\label{equ:compact1}
\dX(x,y) \le\| f(x)-f(y)\|\le C\sqrt{|\log \dX(x,y)|}\dX(x,y).
\end{equation}

The original net-argument from \cite[Proposition 7.18]{BenyaminiLindenstrauss00} induces the loss of a logarithmic factor, but with a bit more care it can be improved. The proof that we give below uses a slightly different embedding.    

\begin{prop}\label{prop:compact2} Let $\metXd$ be a compact metric space. For every continuous, decreasing function $\mu\colon (0,\diam(\metX)]\to [0,+\infty)$ such that $\mu(\diam(\metX))=0$ and $\ds\lim_{t\to0}\mu(t)=+\infty$, there exists an embedding $f\colon \metX\to\xbl$ such that for all $x,y\in \metX$,
\begin{equation*}
\frac{1}{2}\frac{\dX(x,y)}{\mu(\dX(x,y))} \le\| f(x)-f(y)\|\le \frac{\pi}{\sqrt{6}}\dX(x,y).
\end{equation*}
\end{prop}

\begin{proof}
Assume that the diameter of $\metX$ is $D$. And let $\mu\colon (0,D]\to [0,+\infty)$ be a continuous, decreasing function such that $\mu(D)=0$ and $\ds\lim_{t\to0}\mu(t)=+\infty$. Then the map $\mu$ has an inverse denoted $\mu^{-1}\colon [0,+\infty)\to (0,D]$ that is decreasing, with $\ds\lim_{t\to+\infty}\mu^{-1}(t)=0$. Fix $x_0$ in $\metX$, and denote $\sigma:=\mu^{-1}\colon \bZ^-\to (0,D]$. For any $k\in\bZ^+$, let $R_k$ be a maximal $\frac{\sigma(k)}{16}$-net of $\metX$. Observe that $R_k$ is finite by compactness of $\metX$.  Define the following $1$-Lipschitz maps:

$$\begin{array}{rcl}
\varphi_k \colon \metX & \rightarrow & \ell_\infty(R_k)\\
       &   &  \\
     x & \mapsto & \big(\dX(x,y)-\dX(y,x_0)\big)_{y\in R_k}. 
    \end{array}$$
The embedding is given by
$$\begin{array}{rcl}
     f\colon \metX  & \rightarrow &\displaystyle \left(\sum_{k=1}^{\infty}\ell_\infty(R_k)\right)_2\\
       &   &  \\
     x & \mapsto &\displaystyle \sum_{k=1}^\infty\frac{\varphi_k(x)}{k}.\\
    \end{array}$$
It is clear that $f$ is $C$-Lipschitz with $C=(\sum_{k=1}^\infty\frac{1}{k^2})^{\frac12}=\frac{\pi}{\sqrt{6}}$.\\

\smallskip 

\noindent Let $x\neq y\in \metX$, then there exists $l\in\bZ^+$ such that $\sigma(l+1)\le \dX(x,y)<\sigma(l)$, or equivalently $l+1\le \mu(\dX(x,y))< l$. Since $R_{l+1}$ is a $\frac{\sigma(l+1)}{4}$-net in $\metX$ we can find $r_{y}\in R_{l+1}$ such that $\dX(r_{y},y)<\frac{\sigma(l+1)}{4}\le \frac{\dX(x,y)}{4}$. Therefore
\begin{align*}
   \| f(x)-f(y)\| & \ge \frac{\|\varphi_{l+1}(x)-\varphi_{l+1}(y)\|_{\infty}}{l+1} \ge \frac{\sup_{r\in R_{l+1}}\big|\dX(x,r)-\dX(y,r)\big|}{l+1}\\
   &\ge \frac{|\dX(x,r_{y})-\dX(y,r_{y})|}{l+1}\ge \frac{\dX(x,y)-2\dX(y,r_{y})}{l+1}\ge \frac{1}{2}\frac{\dX(x,y)}{l+1},\\
  \end{align*}
and
\begin{equation*}\label{equ3}
\| f(x)-f(y)\|\ge \frac{1}{2}\frac{\dX(x,y)}{\mu(\dX(x,y))}.
\end{equation*}
\end{proof}

Proposition \ref{prop:compact2} says that we can construct uniform embeddings which are arbitrarily close to bi-Lipschitz embeddings. Proposition \ref{prop:compact2} is optimal since there exists a compact metric space which does not admit a bi-Lipschitz embedding into any Banach space with the Radon-Nikod\'ym property, and in particular into any reflexive Banach space (c.f. \cite{BaudierLancien15} for more details).

In \cite{Baudier12}, the original net-argument from \cite[Proposition 7.18]{BenyaminiLindenstrauss00} was combined with barycentric gluing to show that every proper metric space admits a $[t/\log^2 t,t]$-embedding into any Banach space with trivial cotype. This embedding which is simultaneously uniform and coarse was quite satisfactory but was suboptimal in two ways. On the one hand, it is not a coarse bi-Lipschitz embedding, but we know that every proper space admits a coarse bi-Lipschitz embedding into any Banach space which contains the $\ell_\infty^n$'s. On the other hand, the compact balls only embed with suboptimal faithfulness in light of either \eqref{equ:compact1} or Proposition \ref{prop:compact2}. This suboptimality was due to the fact that the parameters of the nets and the radii of the dyadic balls used in the barycentric gluing were correlated in \cite{Baudier12}. By ``decorrelating" the construction, a much tighter embedding result was proved in \cite{BaudierLancien15} and the following notion naturally arose. 

\begin{defi}\label{def:almostLip} 
We say that $(X,\dX)$ is \emph{almost bi-Lipschitzly embeddable} into $\metYd$ if there exist a scaling factor $r\in(0,\infty)$ and a constant $D\in[1,\infty)$, such that for any continuous function $\varphi\colon [0,+\infty)\to [0,1)$ satisfying $\varphi(0)=0$ and $\varphi(t)>0$ for all $t>0$, there exists a map $f_\varphi\colon \metX\to \metY$ such that for all $x,y\in \metX$,
$$\varphi(\dX(x,y))r\dX(x,y)\le \dY(f_\varphi(x),f_\varphi(y))\le Dr\dX(x,y).$$
\end{defi}

Proposition \ref{prop:compact2} essentially says that a compact metric space almost bi-Lipschiztly embeds into $\xbl$. It is clear from Definition \ref{def:almostLip} that if $\metX$ admits a bi-Lipschitz embedding into $\metY$, then $\metX$ almost bi-Lipschitzly embeds into $\metY$. Note also that if $\metX$ almost bi-Lipschitzly embeds into $\metY$, then $\metX$ admits an embedding into $\metY$ that is simultaneously coarse and uniform.  Moreover, if $\metX$ almost bi-Lipschitzly embeds into $\metY$, it is easy to see, by taking $\varphi$ appropriately, that $\metX$ admits a coarse bi-Lipschitz embedding into $\metY$. The following theorem is the main result from \cite{BaudierLancien15}.

\begin{theo}\label{thm:almost}
Let $p\in[1,\infty]$. If $\ell_p$ is finitely representable into a Banach space $\banY$ then any proper subset of $L_p[0,1]$ is almost bi-Lipschitzly embeddable into $\banY$. In particular, any proper subset of $L_p[0,1]$ is almost bi-Lipschitzly embeddable into $\ell_p$.
\end{theo}

The following corollaries can be easily derived from Theorem \ref{thm:almost}.

\begin{coro}\label{cor:proper-cotype}
Let $\banY$ be a Banach space which contains the $\ell_\infty^n$'s, and $\metX$ be a proper metric space, then $\metX$ almost bi-Lipschitzly embeds into $\banY$. In particular, every proper metric space almost bi-Lipschitzly embeds  into the reflexive and asymptotically-$\ell_2$ Banach space \xbl.
\end{coro}

\begin{coro}\label{cor:proper-H}
Let $\banY$ be an infinite dimensional Banach space. Then any proper subset of Hilbert space almost bi-Lipschitzly embeds into $\banY$.
\end{coro}

Theorem \ref{thm:almost}, Corollary \ref{cor:proper-cotype}, Corollary \ref{cor:proper-H} are tight in the sense that ``almost bi-Lipschitz embeddability'' cannot be upgraded to ``bi-Lipschitz embeddability'' (c.f. \cite{BaudierLancien15} for more details).
 
\subsection{Nearly isometric embeddability of stable metric spaces}\label{sec:stablemetric}
Comparing Corollary \ref{cor:BL} and Corollary \ref{cor:proper-cotype} we see that it is always possible to embed locally finite into any Banach space that does not contain the $\ell_\infty^n$'s, and also proper metrics albeit with a lesser degree of faithfulness. It is natural to wonder whether this trend continues if we enlarge once more the class of metrics considered. In order to relax the properness condition we can consider stable metrics. Before discussing embedding results for stable metrics, we make a little detour to introduce the fascinating notion of metric stability.

\subsubsection{Stable metrics}

Stable norms\footnote{A norm is stable if for any two bounded sequences $\xn$, $\yn$, $\lim_{m}\lim_{n}\|x_m+y_n\|=\lim_{n}\lim_{m}\|x_m+y_n\|$, whenever the limits exist.} were introduced, originally for separable Banach spaces, by Krivine and Maurey \cite{KrivineMaurey81}. Whether every infinite-dimensional Banach space contains an isomorphic copy of $\co$ or $\ell_p$ for some $p\in[1,\infty)$ was a long standing open problem in Banach space theory. In early 1970's, Tsirelson \cite{Tsirelson74} built an example of a space that does not have this property. However, Krivine and Maurey showed that for stable Banach spaces such Tsirelson-type construction cannot happen since every stable Banach space contains an isomorphic copy (actually almost isometric) of $\ell_p$ for some $p\in[1,\infty)$. Krivine-Maurey result is one of many evidences that stable Banach spaces have a much regular structure than arbitrary Banach spaces. We point the reader to \cite{Guerre_book92} for an extensive account of the theory of stable Banach spaces. It seems that the natural extension of the notion of stability to arbitrary metric spaces was first studied by Garling \cite{Garling82}.
\begin{defi} A metric space $(X,\dX)$ is said to be \emph{stable} if for any two bounded sequences $\xn$, $\yn$, and any two non-principal ultrafilters $\cU,\cV$ on $\bN$, the following equality holds
\begin{equation*}\label{eq:stability}
\lim_{m\to\cU}\lim_{n\in\cV}d_{X}(x_m,y_n)=\lim_{n\in\cV}\lim_{m\to\cU}d_{X}(x_m,y_n).
\end{equation*} 
\end{defi}
A Banach space is stable if its canonical metric induced by its norm is stable. Stability is an isometric property and is inherited by subsets. Despite condition \eqref{eq:stability} seems quite restrictive the class of stable metrics is rather rich and contains lots of interesting metrics. 
 
\begin{exam}\label{exam:trivial_stable}
Finite, compact, bounded geometry, locally finite, or proper metric spaces are stable.
\end{exam}
For all the metric spaces above the closed balls are either finite or compact, and if $\xn$, $\yn$ are bounded sequences, then given any non-principal ultrafilters $\cU,\cV$ on $\bN$, there exist $x,y\in \metX$ such that $\lim_{n\to\cU}x_n=x$ and $\lim_{n\in\cV}y_n=y$. By continuity properties of the distance function we thus have, 
\begin{equation*}
\lim_{m\to\cU}\lim_{n\in\cV}d_{\metX}(x_m,y_n)=d_{\metX}(x,y)=\lim_{n\in\cV}\lim_{m\to\cU}d_{\metX}(x_m,y_n).
\end{equation*}
In particular, finite-dimensional Banach spaces, finitely generated groups equipped with their canonical word metric, compactly generated groups equipped with their canonical proper metric are stable metric spaces. It is a classical fact that if a Banach space is proper then it must be finite-dimensional. We now describe infinite-dimensional Banach spaces that are stable since it will provide us with examples of stable metric spaces which do not belong to the list of (trivially) stable spaces from Example \ref{exam:trivial_stable}.
\begin{exam} 
Hilbert space is stable.
\end{exam}
It is fairly easy to show that Hilbert space is stable using the classical representation of Hilbert norm in terms of the scalar product. Assume that $\xn$, $\yn$ are two bounded sequences in Hilbert space. Since Hilbert space is reflexive the two sequences are weakly-convergent to, say, $x$ and $y$, respectively. Denote $a=\lim_{n\to\cU}\|x_n\|_2^2$ and $b=\lim_{n\to\cU}\|y_n\|_2^2$. Then for any two ultrafilters $\cU$ and $\cV$, 
\begin{align*}
\lim_{n\to\cU}\lim_{m\in\cV}\|x_n-y_m\|_2^2&=\lim_{n\to\cU}\lim_{m\in\cV}\left(\|x_n\|_2^2+\|y_m\|_2^2-2\langle x_n,y_m\rangle\right)\\
&=\lim_{n\to\cU}\|x_n\|_2^2+\lim_{m\in\cV}\|y_m\|_2^2-2\langle x,y\rangle\\
&=a+b-2\langle x,y\rangle,
\end{align*}
where in the second equality we used the definition of weak-convergence.
Proceeding in a similar manner we can show that
\begin{align*}
\lim_{m\in\cV}\lim_{n\to\cU}\|x_n-y_m\|_H^2=a+b-2<x,y>,
\end{align*} 
and this completes the proof.

\begin{exam}\label{exa:lp-stable}
The Banach space $\ell_p$ for $p\in[1,\infty)$ is stable.
\end{exam}

There are several different approaches to prove the statement in Example \ref{exa:lp-stable}. It is a elementary, but somewhat tedious, exercise to show that for $p\in[1,\infty)$, the $\ell_p$-sum of stable Banach spaces is stable, and one can then argue that $\ell_p$ is by definition the $\ell_p$-sum of countably many copies of $(\bR,|\cdot|)$, which is a proper, and hence stable, space. This argument builds on the following classical lemma. 
\begin{lemm}\label{lem:p-asystu}
Let $\xn$ be a bounded sequence in $\ell_p$, and $\cU$ a non-principal
ultrafilter on $\bN$. Suppose that $\lim_{n\to\cU}e_i^*(x_n)=0$ for all $i\in\bN$, i.e. $\xn$ converges coordinatewise to $0$ with respect to $\cU$. Then, for every $z\in\ell_p$,
\begin{equation}\label{eq:p-asystr}
\lim_{n\to\cU}\|z+x_n\|_p^p=\|z\|_p^p+\lim_{n\to\cU}\|x_n\|_p^p.
\end{equation}
\end{lemm}
The conclusion of the lemma clearly holds if $x$ and the $x_n$'s have disjoint supports. Since $\xn$ converges coordinate-wise to $0$ by reaching far out in the sequence one can select $x_n$ such that the essential contributions to the norm of $x_n$ and of $x$ are supported on essentially disjoint supports, and a classical approximation and truncation argument gives the conclusion. 

For every bounded sequence $\yn$ and every non-principal ultrafilter $\cU$ on $\bN$, there exist $y\in \ell_p$ and $\mu\ge 0$ so that if $z\in\ell_p$, then $\lim_{n\to\cU}\|z-y_n\|_p^p=\|z-y\|^p+\mu^p$. Indeed, there exists $y\in\ell_p$ such that for all $i\in\bN$ $\lim_{n\to\cU}e_i^*(y_n)=e_i^*(y)$ and the conclusion follows from Lemma \ref{lem:p-asystu} with $\xn=(y-y_n)_{n=1}^\infty$, and $\mu=\lim_{n,\cU}\|y-y_n\|$. We are now in position to prove that $\ell_p$ is stable. For a pair of bounded sequences $\xn$, $\yn$, and a pair of non-principal ultrafilters $\cU,\cV$ on $\bN$, there exist $x,y\in \ell_p$, $\mu,\nu\ge 0$ such that for all $z\in \ell_p$ the following equalities hold: 
$$\lim_{n\to\cU}\|z-x_n\|_p^p=\|z-x\|_p^p+\mu^p$$
and
$$\lim_{n\in\cV}\|z-y_n\|_p^p=\|z-y\|_p^p+\nu^p.$$
It follows that $$\lim_{m\in\cV}\lim_{n\to\cU}\|x_n-y_m\|_p^p=\|x-y\|_p^p+\mu^p+\nu^p=\lim_{n\to\cU}\lim_{m\in\cV}\|x_n-y_m\|_p^p.$$

\begin{exam}\label{exam:co}
The Banach space $\co$, and more generally any Banach space containing an isomorphic copy of $\co$, is not stable.
\end{exam}

The fact that $\co$ is not stable can easily be checked directly by considering the sequences $x_n=-e_n$ and $y_n=\sum_{i=1}^ne_i$, $n\geq 1$ where $\en$ denotes the canonical basis of $\co$. Indeed, $\lim_{m\in\cV}\lim_{n\to\cU}\|x_n-y_m\|_\infty=2$ while $\lim_{n\to\cU}\lim_{m\in\cV}\|x_n-y_m\|_\infty=1$. In order to show the second part of Example \ref{exam:co} we need to invoke James' $\co$-distortion theorem which says that if $Y$ contains an isomorphic copy of $\co$ then $Y$ will contain for every $\epsilon>0$ a subspace that is $(1+\epsilon)$-isomorphic to $\co$. More examples of spaces that are not stable will be given in Section \ref{sec:upper-stability}.

The stability of the function space $L_p[0,1]$ is much more difficult to obtain than the stability of the sequence space $\ell_p$ since the validity of equality \eqref{eq:p-asystr} for the sequence space does not hold in the function space case, and thus another argument is needed.

\begin{exam}\label{exam:L_p}
The Banach space $L_p[0,1]$ for $p\in[1,\infty)$ is stable.
\end{exam}

For $s\in(0,1)$, the metric space $(\metX,\sd^s_{\metX})$ is commonly called the $s$-\emph{snowflaking} of $\metX$. It is clear that a snowflaking of a metric pace is stable if and only if the original metric on the space is stable. For $p\in[1,2)$, $L_p[0,1]$ will be stable since it is well known that the $\frac{p}{2}$-snowflaking of $L_p[0,1]$ embeds isometrically into the Hilbert space $L_2[0,1]$. This argument fails for $p>2$, and the stability of $L_p[0,1]$ for $p>2$ is more difficult to prove and requires a fine understanding of the relationship between stability and reflexivity. Krivine and Maurey gave a representation of the norm of a stable Banach space which provides a direct relationship between stability and reflexivity\footnote{The connection between stable functions and weak-compactness goes back at least to the work of Grothendieck \cite{Grothendieck52}.}. 

\begin{theo}\label{represent}
Let $\banX$ be a Banach space and fix $p\in[1,\infty)$. Then, $\banX$ is stable if and only if there exist a reflexive Banach space $\banY$, a dense subset $B$ of the unit ball of $\banX$, and maps $g\colon B\to \banY$, $h\colon B\to \banY^*$ so that for all $x,y
\in B$ we have $\|x-y\|^p=\langle g(x),h(y)\rangle$ where $\langle \cdot,\cdot\rangle$ is the duality product between $\banY$ and $\banY^*$.
\end{theo}
The theorem below, which includes Example \ref{exam:L_p}, follows from a similar but more delicate representation of stable norms that was also proved in \cite{KrivineMaurey81}. For a Banach space $\banX$ we denote by $L_p(\Omega,\cB, \mu;\banX)$, or simply $L_p(\Omega;\banX)$,  the Banach space of Bochner equivalence classes of Bochner $p$-integrable, $\banX$-valued functions defined on the measured space $(\Omega, \cB, \mu)$. 

\begin{theo}
Let $p\in[1,\infty)$. If $\banX$ is a stable Banach space then $L_p([0,1];\banX)$ is also stable.
\end{theo} 

\subsubsection{Kalton's embedding of stable metrics}
It follows from the examples in the previous section that the class of stable metrics strictly extends the class of proper metrics. For stable metrics there is no analogue of Corollary \ref{cor:BL} or Corollary \ref{cor:proper-cotype} even if we are willing to significantly weaken the degree of faithfulness of the embeddings. Indeed, it was proved in \cite{BLS_JAMS} that Tsirelson space $T^*$ (which is reflexive and contains the $\ell_\infty^n$'s) does not contain a coarse copy nor a uniform copy of any of the (stable) Banach spaces $\ell_p$ for $p\in[1,\infty)$. The next course of action would be either to consider even weaker notions of embeddability or to settle for embeddings into Banach spaces with weaker geometric features.

A remarkable result of Kalton \cite{Kalton07} states that a stable metric can always be embedded into a reflexive Banach space, which depends on the stable metric, in a way that distorts the distances by only a slight amount.

\begin{theo}\label{thm:Kalton}
Let $\metX$ be a stable metric space, and $s\in(0,1)$, then there exist a reflexive Banach space $\ban K_s(\metX)$ and a map $f\colon \metX\to \ban K_s(\metX)$ such that for all $x,y\in \metX$
\begin{equation*}
\min\{\dX(x,y),\sd^s_\metX(x,y)\}\le \|f(x)-f(y)\|_{\ban K}\le\max\{\dX(x,y),\sd^s_\metX(x,y)\}.
\end{equation*}
\end{theo}

Unlike in Corollary \ref{cor:BL} or Corollary \ref{cor:proper-cotype} where the host space can be taken arbitrarily from a given class of Banach spaces, the reflexive space $\ban K_s(\metX)$, and the embedding $f$, highly depend on $s$ and $\metX$. Given $s\in(0,1)$ the crux of Kalton's proof is to leverage the stability condition in order to construct a weakly compact operator $S\colon \ell_1(W)\to \lip_0(\metX,\sd^s_\metX)$ for some well-chosen set $W$ that is weakly relatively compact in $\lip_0(\metX,\sd^s_\metX)$ the Lipschitz space over the $s$-snowflaking of $\metX$. Then, by Davis-Figiel-Johnson-Pelczy\'nski factorization theorem \cite{DFJP74}, $S$ factors through a reflexive Banach space $\banX$ and $\ban K_s(\metX)$ is defined as the dual of $\banX$. A careful inspection of Kalton's proof reveals that an embedding with tighter guarantees can be achieved. The adjustments were carried over in \cite{BaudierLancien15} were the following definition, which captures the quantitative improvements, was introduced.

\begin{defi} Let $\cC$ be a class of metric spaces. We say that $\metXd$ \emph{nearly isometrically embeds} into $\cC$ if for any pair of continuous functions $\rho, \omega:[0,+\infty)\to [0,+\infty)$ satisfying
\begin{enumerate}[(i)]
\item $\omega(0)=0$, $t\le \omega(t)$ for $t\in [0,1]$, and $\lim_{t\to 0}\frac{\omega(t)}{t}=+\infty$, 
\item $\omega(t)=t$ for $t\ge 1$,
\item $\rho(t)=t$ for $t\in [0,1]$,
\item $\rho(t)\le t$ for $t\ge 1$ and $\lim_{t\to +\infty}\frac{\rho(t)}{t}=0$,
\end{enumerate}
there exist a space $\metYd\in \cC$ and a map $f\colon \metX\to \metY$ such that for all $x,y\in \metX$
$$ \rho(\dX(x,y))\le \dY(f(x),f(y)) \le \omega(\dX(x,y)).$$
If $\cC$ is reduced to a single element $\metYd$ we say that $\metXd$ nearly isometrically embeds into $\metYd$.
\end{defi}

\begin{theo}\label{nearlyisom}
Let $\metXd$ be a stable metric space, then $\metX$ nearly isometrically embeds into the class of reflexive Banach spaces. 
\end{theo}
\noindent We refer to \cite{BaudierLancien15} for the proof of Theorem \ref{nearlyisom}. The terminology is motivated by the fact that if $\metX$ nearly isometrically embeds into $\cC$, then for every $0<\delta\le \Delta<\infty$ there exist $\metY\in \cC$ and $f\colon \metX\to \metY$ such that 
for all $x,y\in \metX$ satisfying $\dX(x,y)\in[\delta,\Delta]$,
$$\dY(f(x),f(y))=\dX(x,y).$$ This property is not achieved with Kalton's original embedding.

\section{Upper stability}\label{sec:upper-stability}

Krivine-Maurey stability theorem can be used to provide examples of Banach spaces without any stable renormings.

\begin{theo}\label{thm:KM}
Every stable Banach space contains an isomorphic copy\footnote{In fact it contains an almost linearly isometric copy.} of $\ell_p$ for some $p\in[1,\infty)$
\end{theo}

Therefore, Tsirelson-like spaces, which do no contain any isomorphic copy $\ell_p$ for any $p\in[1,\infty)$, do not admit any equivalent stable norm. Note also that if a Banach space contains an isomorphic copy of $\ell_p$ for some $p\in[1,\infty)$, then it will contain an unconditional basic sequence. Hence, Banach spaces with no unconditional basic sequences, in particular hereditarily indecomposable Banach spaces (see \cite{Maurey03} for a discussion of how to construct such spaces), also do not admit any equivalent norm that is stable. 

In \cite{Raynaud83}, Raynaud showed that the conclusion of Krivine-Maurey theorem is valid for any Banach space whose unit ball admits a uniformly equivalent stable metric. It was well known that spreading sequences in stable Banach spaces are unconditional, and Raynaud showed that this property still holds if we merely assume that the unit ball of the Banach space embeds uniformly into a stable space. 

\begin{theo}\label{thm:spreading-unconditional}
Let $\banX$ be a Banach space. If the unit ball of $\banX$ uniformly embeds into stable metric space, then every spreading basic sequence in $\banX$ is unconditional.
\end{theo}

The following theorem for \cite{Raynaud83} follows then from the observation that the summing basis of $\co$, i.e. the sequence $\{s_n\}_n$ where $s_n=\sum_{i=1}^n e_i$, is spreading but not unconditional, an extends a previous result of Enflo \cite{Enflo69} which states that the unit ball of $\co$ does not embed uniformly into Hilbert space.

\begin{theo}\label{thm:co-stable}
The unit ball of $\co$ does not uniformly embed into a stable metric space.
\end{theo}
The notion of stability made a spectacular come back in the nonlinear embedding due to the work of Kalton in \cite{Kalton07} where Property $\cQ$ was introduced. In the sequel we will always write an element $\nbar=(n_1,\dots,n_k)$ of $\Nk$, the set of $k$-subsets of $\bN$, in increasing order, i.e. $n_1<n_2<\cdots<n_k$. In \cite{Kalton07}, Kalton defined a graph structure on (the vertex set) $\Nk$ as follows: $\mbar$ and $\nbar$ are adjacent if and only if they interlace, meaning that either $m_1\le n_1\le m_2\le n_2\le\cdots\le m_k\le n_k$ or  $n_1\le m_1\le n_2\le m_2\le\cdots\le n_k\le n_k$. This (non-locally finite) infinite graph is then equipped with its canonical graph metric, simply denoted $\sd_{\sI}$. Kalton's property $\cQ$ is the following concentration inequality on Kalton's interlacing graphs $\{\Nk,\sd_{\sI}\}_{k\ge 1}$.

\begin{defi} A metric space $\metXd$ has property $\cQ$ if there exists a universal constant $C>0$ so that for every Lipschitz map $f\colon (\Nk,\sd_\sI)\to \metX$ there is an infinite subset $\bM$ of $\bN$, such that for all $\mbar,\nbar\in \Mk{\bM}{k}$,
\begin{equation}\label{eq:propertyQ}
\dX(f(\mbar),f(\nbar))\le C\cdot\lip(f).
\end{equation}
\end{defi}

Note that for any infinite subset $\bM$ of $\bN$, the graph $(\Mk,\sd_\sI)$ has diameter $k$, and it is immediate from the definitions that a metric space with property $\cQ$ cannot equi-coarsely  contain the sequence $\{(\Nk,\sd_\sI)\}_{k\ge 1}$. Kalton's paper \cite{Kalton07} was extremely influential. In particular, it initiated the use of concentration inequalities in the form of \eqref{eq:propertyQ} in nonlinear embedding problems. In \cite{KaltonLova08}, \cite{Kalton11}, or \cite{BLS_JAMS} several long-standing open problems were resolved using this approach.

Kalton's interlacing graph metric encapsulates the particular behavior of the summing basis of $\co$ as it is a combinatorial realization of a metric that is naturally induced by the bi-monotone version of the summing norm given by 
\begin{align}
\label{summing norm}
\left\|\sum_ia_is_i\right\|_{\mathrm{sum}}= \sup\left\{\left|\sum_{i=k}^ma_i\right|: k,m\in \bN,\,  k\le m\right\}.
\end{align}

It is a tedious task (cf. \cite{LPP18} or \cite{BLMS_FM}) to show that if $\mbar,\nbar\in\Nk$ then 
\begin{equation}
\sd_{\sI}(\nbar,\mbar)=\left\|\sum_{i=1}^k s_{m_i}-\sum_{i=1}^k s_{n_i}\right\|_{sum}.
\end{equation}

It was shown in \cite{Kalton07} that Kalton's property $\cQ$ belongs to the scanty list of obstructions to coarse and uniform embeddability.

\begin{lemm}\label{lem:Q}
Let $\banX$ be a Banach space. 
Assume that, either
\begin{enumerate}[(i)]
\item $\banX$ coarsely embeds into a metric space with property $\cQ$
\item[] or
\item the unit ball of $\banX$ uniformly embeds into a metric space with property $\cQ$.
\end{enumerate}
Then $\banX$ has property $\cQ$.
\end{lemm}


Since Kalton showed that stable metric spaces have property $\cQ$, Raynaud's non-embeddability result (and its coarse version) follows from the discussion above which implies that $\co$ does not have property $\cQ$, and Lemma \ref{lem:Q}. The arguments of Raynaud and Kalton rely on the following crucial observation (we refer to \cite{Raynaud83} for a proof).

\begin{lemm}Let $\metXd$ be a stable metric space. Let $1\le \ell \le k$, $\pi\colon [k]\to [k]$ a permutation which preserves the order on $\{1,\dots,\ell\}$ and $\{\ell+1,\dots,k\}$, and $\cU_1,\dots,\cU_k$ be ultrafilters on $\bN$. Then for all bounded maps $f\colon \Sk{\bN}{\ell}\to \metX$ and $g\colon \Sk{\bN}{k-\ell}\to \metX$ the following identity holds
\begin{eqnarray*}
\lim_{n_1\to\cU_1}\dots\lim_{n_k\to\cU_k}\dX(f(n_1,\dots,n_\ell),g(n_{\ell+1},\dots,n_k))=\\
\lim_{n_1\to\cU_1}\dots\lim_{n_k\to\cU_k}\dX(f(n_{\pi(1)}, \dots,n_{\pi(\ell)}),g(n_{\pi(\ell+1)},\dots,n_{\pi(k)})
\end{eqnarray*}
\end{lemm}

The notion of upper stability that we introduce in the definition below is a natural bi-Lipschitz variant of the isometric notion of stability and can be used to quantify the lack of stability. 

\begin{defi}
Let $\metXd$ be a metric space. We say that $\metX$ is \emph{$K_u$-upper stable} if  for every $k\ge 1$, every $1\le \ell \le k$, every permutation $\pi\colon [k]\to [k]$ which preserves the order on $\{1,\dots,\ell\}$ and $\{\ell+1,\dots,k\}$, and every bounded maps $f\colon \Sk{\bN}{\ell}\to \metX$ and $g\colon \Sk{\bN}{k-\ell}\to \metX$ we have for every infinite subset $\bM$ of $\bN$ 
\begin{equation}\label{eq:upper-stable}
\inf_{\nbar\in[\bM]^k}\dX(f(n_1,\dots,n_\ell),g(n_{\ell+1},\dots,n_k))\le K_u \sup_{\nbar\in [\bM]^k}\dX(f(n_{\pi(1)},\dots,n_{\pi(\ell)}),g(n_{\pi(\ell+1)},\dots,n_{\pi(k)}).
\end{equation}
And we say that $\metX$ is \emph{upper stable} if it is $K_u$-upper stable for some $K_u>0$.  
\end{defi}

The bi-Lipschitz nature of upper stability makes it a more convenient and natural notion to work with since results for stable metrics that are of isomorphic nature can be extended to the upper stable metrics. The first observation is that upper stability generalizes stability.

\begin{lemm}
Every stable metric space is 1-upper stable.
\end{lemm}

\begin{proof}
Let $\bM$ be an infinite subset of $\bN$ and $\cU_1,\dots,\cU_k$ be ultrafilters on $\bN$ containing $\bM$.
Since for every $\nbar\in[\bM]^k$,

$$\dX(f(n_{\pi(1)},\dots,n_{\pi(\ell)}),g(n_{\pi(\ell+1)},\dots,n_{\pi(k)})\le \sup_{\nbar\in [\bM]^k}\dX(f(n_{\pi(1)},\dots,n_{\pi(\ell)}),g(n_{\pi(\ell+1)},\dots,n_{\pi(k)})$$

and since the ultrafilters contain $\bM$ we have

\begin{eqnarray*}
\lim_{n_1\to\cU_1}\dots\lim_{n_k\to\cU_k}\dX(f(n_{\pi(1)},\dots,n_{\pi(\ell)}),g(n_{\pi(\ell+1)},\dots,n_{\pi(k)}) \le&\\
\sup_{\nbar\in [\bM]^k}\dX(f(n_{\pi(1)},\dots,n_{\pi(\ell)}),g(n_{\pi(\ell+1)},\dots,n_{\pi(k)}).&
\end{eqnarray*}

By stability,

\begin{eqnarray*}
\lim_{n_1\to\cU_1}\dots\lim_{n_k\to\cU_k}\dX(f(n_1,\dots,n_\ell),g(n_{\ell+1},\dots,n_k))\le &\\
\sup_{\nbar\in [\bM]^k}\dX(f(n_{\pi(1)},\dots,n_{\pi(\ell)}),g(n_{\pi(\ell+1)},\dots,n_{\pi(k)}),&
\end{eqnarray*}

and thus for every $\vep>0$ there exists $\nbar\in [\bM]^k$ (we again use that the ultrafilters contain $\bM$) such that 

$$\dX(f(n_1,\dots,n_\ell),g(n_{\ell+1},\dots,n_k))\le \sup_{\nbar\in [\bM]^k}\dX(f(n_{\pi(1)},\dots,n_{\pi(\ell)}),g(n_{\pi(\ell+1)},\dots,n_{\pi(k)})+\vep$$

and the conclusion follows.
\end{proof}

Akin Property $\cQ$, upper stability is clearly a bi-Lipschitz invariant of metric spaces, but also a coarse and uniform invariant for Banach spaces.

\begin{lemm}\label{lem:coarse-upper}
If a Banach space $\banX$ coarsely, or uniformly, embeds into an upper stable metric space, then $\banX$ is upper stable.
\end{lemm}

\begin{proof}
We prove first the coarse statement.
Since $\banXn$ coarsely embeds into a stable metric space $\metYd$ there exists $\varphi\colon \banX\to \metY$, and $\rho,\omega\colon [0,\infty)\to [0,\infty)$ such that 
\begin{equation*}
\norm{x-y}\le \dY(\varphi(x),\varphi(y))\le\omega(\norm{x-y}).
\end{equation*}

Let $1\le \ell \le k$, $f\colon \Sk{\bN}{\ell}\to \banX$ and $g\colon \Sk{\bN}{k-\ell}\to \banX$ bounded maps, and $\pi\colon [k]\to [k]$ a permutation that preserves the order on $\{1,\dots,\ell\}$ and $\{\ell+1,\dots,k\}$, be given. Let $\bM$ and infinite subset of $\bN$ and let $$\alpha\eqd\sup_{\nbar\in\Mk} \norm{f(n_{\pi(1)},\dots,n_{\pi(\ell)})-g(n_{\pi(\ell+1)},\dots,n_{\pi(k)})}.$$
Then for every $(n_1,\dots, n_k)\in \Mk$

\begin{eqnarray*}
\dY(\varphi(\frac{1}{\alpha}f(n_{\pi(1)},\dots,n_{\pi(\ell)})),\frac{1}{\alpha}\varphi(g(n_{\pi(\ell+1)},\dots,n_{\pi(k)}))\le\omega(1),
\end{eqnarray*}

and since hence $\omega(1)$ is independent of $n_1<\dots<n_k$ we have

\begin{equation*}
\sup_{\nbar\in \Mk}\dY(\frac{1}{\alpha}\varphi(f(n_{\pi(1)},\dots,n_{\pi(\ell)})),\varphi(\frac{1}{\alpha}g(n_{\pi(\ell+1)},\dots,n_{\pi(k)}))\le \omega(1)
\end{equation*}

Let $\beta\eqd \frac{1}{\alpha}\inf_{\nbar\in\Mk} \norm{f(n_1,\dots,n_\ell)-g(n_{\ell+1},\dots,n_k)}$ and observe that for all $(n_1,\dots,n_k)\in\Mk$

\begin{equation*}
\dY(\varphi(\frac{1}{\alpha}f(n_1,\dots,n_\ell)),\varphi(\frac{1}{\alpha}g(n_{\ell+1},\dots,n_k)))\ge \inf\{\dY(\varphi(x),\varphi(y))\colon \dX(x,y)\ge \beta\}
\end{equation*}

and since the right hand side is independent of $n_1<\dots<n_k$,
 
\begin{equation*}
\inf_{\nbar\in\Mk}\dY(\varphi(\frac{1}{\alpha}f(n_1,\dots,n_\ell)),\varphi(\frac{1}{\alpha}g(n_{\ell+1},\dots,n_k)))\ge \inf\{\dY(\varphi(x),\varphi(y))\colon \dX(x,y)\ge \beta\}
\end{equation*}

The upper stability assumption gives

\begin{equation}\label{eq:contradiction-coarse}
\inf\{\dY(\varphi(x),\varphi(y))\colon \dX(x,y)\ge \beta\}\le K_u\omega(1).
\end{equation}

Let $C>0$ such that $\rho(C)\ge 2 K_u\omega(1)$ (such an $C$ exists since by assumption $\lim_{t\to\infty}\rho(t)=\infty$). If  
$\beta > C$ then whenever $\dX(x,y)\ge \beta$ we have $\dY(\varphi(x),\varphi(y))\ge \rho(\dX(x,y))\ge \rho(C)\ge 2K_u\omega(1)$, but this contradicts \eqref{eq:contradiction-coarse}. Therefore, $\beta\le C$ necessarily and this completes the proof since it precisely means that 
\begin{equation*}
\inf_{\nbar\in[\bM]^k}\dX(f(n_1,\dots,n_\ell),g(n_{\ell+1},\dots,n_k))\le C \sup_{\nbar\in [\bM]^k}\dX(f(n_{\pi(1)},\dots,n_{\pi(\ell)}),g(n_{\pi(\ell+1)},\dots,n_{\pi(k)}).
\end{equation*}

Now the uniform case. Let $\beta\eqd \inf_{\nbar\in\Mk} \norm{f(n_1,\dots,n_\ell)-g(n_{\ell+1},\dots,n_k)}$. Then for every $(n_1,\dots, n_k)\in \Mk$

\begin{equation*}
\dY(\varphi(\frac{1}{\beta}f(n_1,\dots,n_\ell)),\varphi(\frac{1}{\beta}g(n_{\ell+1},\dots,n_k)))\ge \rho(1)
\end{equation*}

and since $\rho(1)$ is independent of $n_1<\dots<n_k$ we have

\begin{equation*}
\inf_{n\in\Mk}\dY(\varphi(\frac{1}{\beta}f(n_1,\dots,n_\ell)),\varphi(\frac{1}{\beta}g(n_{\ell+1},\dots,n_k)))\ge \rho(1)
\end{equation*}

Let $$\alpha\eqd\frac{1}{\beta}\sup_{\nbar\in\Mk} \norm{f(n_{\pi(1)},\dots,n_{\pi(\ell)})-g(n_{\pi(\ell+1)},\dots,n_{\pi(k)})},$$

and observe that for all $(n_1,\dots,n_k)\in\Mk$

\begin{equation*}
\dY(\varphi(\frac{1}{\beta}f(n_{\pi(1)},\dots,n_{\pi(\ell)}),\varphi(\frac{1}{\beta}g(n_{\pi(\ell+1)},\dots,n_{\pi(k)})))\le\sup\{\dY(\varphi(x),\varphi(y))\colon \dX(x,y)\le\alpha\}
\end{equation*}

and since the right-hand side is independent of $n_1<\dots<n_k$,
 
\begin{equation*}
\sup_{n\in\Mk}\dY(\varphi(\frac{1}{\beta}f(n_{\pi(1)},\dots,n_{\pi(\ell)}),\varphi(\frac{1}{\beta}g(n_{\pi(\ell+1)},\dots,n_{\pi(k)})))\le\sup\{\dY(\varphi(x),\varphi(y))\colon \dX(x,y)\le\alpha\}
\end{equation*}

The upper stability assumption gives

\begin{equation}\label{eq:contradiction-uniform}
\rho(1)\le K_u\sup\{\dY(\varphi(x),\varphi(y))\colon \dX(x,y)\le\alpha\}
\end{equation}

Let $C>0$ such that $\omega(C)\le \frac{\rho(1)}{2 K_u}$ (such an $C$ exists since by assumption $\lim_{t\to 0}\omega(t)=0$). If  
$\alpha < C$ then whenever $\dX(x,y)\le \alpha$ we have $\dY(\varphi(x),\varphi(y))\le \omega(\dX(x,y))\le \omega(C)\le \frac{\rho(1)}{2 K_u}$, but this contradicts \eqref{eq:contradiction-uniform}. Therefore, $\alpha\ge C$ necessarily and this completes the proof since it precisely means that 
\begin{equation*}
\inf_{\nbar\in[\bM]^k}\dX(f(n_1,\dots,n_\ell),g(n_{\ell+1},\dots,n_k))\le \frac{1}{C} \sup_{\nbar\in [\bM]^k}\dX(f(n_{\pi(1)},\dots,n_{\pi(\ell)}),g(n_{\pi(\ell+1)},\dots,n_{\pi(k)}).
\end{equation*}
\end{proof}

It is not clear whether upper stability is preserved under uniform embeddings of the unit ball. The following lemma, which follows from a careful inspection of the proof in the uniform case of Lemma \ref{lem:coarse-upper}, will be needed to handle uniform embeddings of the unit ball.

\begin{lemm}\label{lem:metric-uniform}
Let $\metXd$ be a metric space that uniformly embeds, with compression function $\rho$ and expansion function $\omega$, into an upper $K_u$-stable metric space. If $1\le \ell \le k$, $f\colon \Sk{\bN}{\ell}\to \metX$ and $g\colon \Sk{\bN}{k-\ell}\to \metX$  are bounded maps, and $\pi\colon [k]\to [k]$ is a permutation that preserves the order on $\{1,\dots,\ell\}$ and $\{\ell+1,\dots,k\}$ then for every infinite subset $\bM$ of $\bN$

\begin{eqnarray*}
\sup_{\nbar\in\Mk} \dX(f(n_{\pi(1)},\dots,n_{\pi(\ell)}),g(n_{\pi(\ell+1)},\dots,n_{\pi(k)}))\ge C,
\end{eqnarray*}

where $C$ is any number satisfying $\omega(C)\le \frac{1}{2K_u}\rho(\alpha)$ with $$\alpha\eqd\inf_{\nbar\in\Mk} \dX(f(n_1,\dots,n_\ell),g(n_{\ell+1},\dots,n_k)) .$$
\end{lemm}

We will now extend the results of Raynaud and Kalton to upper stable spaces. We start with Raynaud's result about unconditionaity of spreading basic sequence in stable spaces. First recall some Banach space theoretic concepts mentioned above. A basic sequence $\xn$ in a Banach space $\banXn$ is \emph{spreading} if it is equivalent to all of its subsequences, in the sense that there exist $A_S, B_s\ge 1$ such that for all $k\ge 1$, $a_1,a_2,\dots,a_k\in\bR$, and $1\le m_1<m_2<\dots<m_k\in \bN$ we have
\begin{equation*}
\frac{1}{B_s}\|\sum_{i=1}^k a_i x_{i}\|\le \|\sum_{i=1}^k a_i x_{n_i}\|\le A_s\|\sum_{i=1}^k a_i x_{i}\|.
\end{equation*}

A basic sequence $\xn$ is said to be \emph{unconditional} if there exists $C_u\ge 1$ such that for all $k\ge 1$, $a_1,a_2,\dots,a_k\in\bR$, and $\vep_1, \vep_2,\dots, \vep_k\in \{-1,+1\}$ we have
\begin{equation*}
\frac{1}{C_u}\|\sum_{n=1}^k a_n x_{n}\|\le \|\sum_{n=1}^k \vep_na_n x_n\|\le C_u\|\sum_{n=1}^k a_n x_{n}\|
\end{equation*}

It is plain that the canonical basis $\{e_n\}_{n\ge 1}$ of $\co$ or $\ell_p$, for $1\le p<\infty$, are both spreading and unconditional\footnote{Sequences that are both spreading and unconditional are called subsymmetric.} (both with constant~1). 

Before we give the proof, observe that the following corollary follows from the fact that the summing basis of $\co$ is spreading but not unconditional.

\begin{theo}\label{thm:s-u-coarse}
Let $\banX$ be a Banach space. Assume that, either
\begin{enumerate}[(i)]
\item $\banX$ coarsely embeds into an upper stable metric space
\item[] or
\item the unit ball of $\banX$ uniformly embeds into an upper stable metric space.
\end{enumerate} 
Then every spreading basic sequence in $\banX$ is unconditional.
\end{theo}

\begin{proof}
We begin by proving $(i)$. Assume that the basic sequence $\xn$ is $C_b$-basic, i.e.
\begin{equation*}
\norm{\sum_{i=1}^j a_i x_{i}\|\le C_b\|\sum_{i=1}^k a_i x_{i}}
\end{equation*} 
for all $1\le j\le k$, and  $a_1,\dots,a_k\in\bR$. Fix $a_1,\dots,a_k\in\bR$ and for any $\eta\in\cube{k}$ let 
$$N(\eta)\eqd\norm{\sum_{i=1}^k \eta_i a_i x_i}.$$

For $\vep=(\vep_1,\dots,\vep_k)\in\cube{k}$, we will  show that there exists a universal constant $\gamma$ such that

\begin{equation}\label{eq:uncond}
N(\1)=\norm{\sum_{i=1}^k a_i x_{i}}\le \gamma \norm{\sum_{i=1}^k \vep_i a_i x_{i}}=\gamma N(\vep).
\end{equation} 

Let $\ell\eqd|\{i\colon \vep_i=+1\}|$ and $\pi$ the permutation on $[k]$ that maps $\{1,\dots,\ell\}$ onto $\{i\colon \vep_i=+1\}$ and $\{\ell+1,\dots,k\}$ onto $\{i\colon \vep_i=-1\}$ while preserving the order on the respective sets. If we define for $(m_1,\dots,m_\ell)\in\Nk$ 
$$f(m_1,\dots,m_\ell)=\sum_{i=1}^\ell a_i x_{m_i}$$

and for $(m_1,\dots,m_{k-\ell})\in\Sk{\bN}{k-\ell}$ 

$$g(m_{1},\dots,m_{k-\ell})= \sum_{i=1}^{k-\ell} a_{\ell+i} x_{m_i},$$

it follows from Lemma \ref{lem:coarse-upper} that  

\begin{equation*}
\inf_{\nbar\in\Nk} \norm{f(n_1,\dots,n_\ell)-g(n_{\ell+1},\dots,n_k)}\le C \sup_{\nbar\in\Nk} \norm{f(n_{\pi(1)},\dots,n_{\pi(\ell)})-g(n_{\pi(\ell+1)},\dots,n_{\pi(k)})},
\end{equation*}

where $C$ is any constant such that $\rho(C)\ge 2K_u\omega(1)$. Observe now that for all $(n_1,\dots,n_k)\in\Nk$

\begin{equation*}
N(\vep)\ge\frac{1}{A_s} \norm{\sum_{i=1}^k \vep_i a_i x_{n_i}}=\frac{1}{A_s} \norm{\sum_{i=1}^\ell a_i x_{n_{\pi(i)}}-\sum_{i=\ell+1}^k a_i x_{n_{\pi(i)}}}.
\end{equation*}

and hence 
\begin{equation*}
\sup_{\nbar\in\Nk} \norm{f(n_{\pi(1)},\dots,n_{\pi(\ell)})-g(n_{\pi(\ell+1)},\dots,n_{\pi(k)})}\le A_sN(\vep).
\end{equation*}

On the other hand, for all $\nbar\in \Nk$

\begin{eqnarray*}
N(\1)\le B_s\norm{\sum_{i=1}^k a_i x_{n_i}} & \le & B_s\norm{-\sum_{i=1}^\ell a_i x_{n_i}+\sum_{i=\ell+1}^k a_i x_{n_i}}+2B_s\norm{\sum_{i=1}^\ell a_i x_{n_i}}\\
		& \le & B_s\norm{\sum_{i=1}^\ell a_i x_{n_i}-\sum_{i=\ell+1}^k a_i x_{n_i}}+2B_sC_b\norm{\sum_{i=1}^\ell a_i x_{n_i}-\sum_{i=\ell+1}^k a_i x_{n_i}}\\
		& = & B_s(2C_b+1)\norm{\sum_{i=1}^\ell a_i x_{n_i}-\sum_{i=\ell+1}^k a_i x_{n_i}},
\end{eqnarray*}

and hence

\begin{equation*}
\inf_{\nbar\in\Nk} \norm{f(n_1,\dots,n_\ell)-g(n_{\ell+1},\dots,n_k)}\ge \frac{N(\1)}{B_s(2C_b+1)}.
\end{equation*}

Consequently $N(\1)\le CA_sB_s(2C_b+1)N(\vep)$ which is \eqref{eq:uncond} with $\gamma=CA_sB_s(2C_b+1)$.

\medskip

The proof of $(ii)$ goes as follows. If $N(\1)\le N(\vep)$ there is nothing to prove and in the sequel we assume that $N(\vep)\le N(\1)$.

Let $\ell\eqd|\{i\colon \vep_i=+1\}|$ and $\pi$ the permutation on $[k]$ that maps $\{1,\dots,\ell\}$ onto $\{i\colon \vep_i=+1\}$ and $\{\ell+1,\dots,k\}$ onto $\{i\colon \vep_i=-1\}$ while preserving the order on the respective sets. Observe now that for all $(n_1,\dots,n_k)\in\Nk$

\begin{equation*}
N(\vep)\ge\frac{1}{A_s} \norm{\sum_{i=1}^k \vep_i a_i x_{n_i}}=\frac{1}{A_s} \norm{\sum_{i=1}^\ell a_i x_{n_{\pi(i)}}-\sum_{i=\ell+1}^k a_i x_{n_{\pi(i)}}}.
\end{equation*}

and

\begin{equation*}
\norm{\sum_{i=1}^\ell a_i x_{n_{\pi(i)}}+\sum_{i=\ell+1}^k a_i x_{n_{\pi(i)}}}=\norm{\sum_{i=1}^k a_i x_{n_i}}\le A_s N(\1)
\end{equation*}

Therefore, $\max \{\norm{\sum_{i=1}^\ell a_i x_{n_{\pi(i)}}}, \norm{\sum_{i=\ell+1}^k a_i x_{n_{\pi(i)}}}\} \le A_s N(\1)$.

Since

\begin{equation*}
\norm{\sum_{i=1}^\ell a_i x_{n_i}} \le A_s\norm{\sum_{i=1}^\ell a_i x_{i}}\le A_sC_b N(\1)
\end{equation*}

and 
\begin{eqnarray*}
\norm{\sum_{i=\ell+1}^k a_i x_{n_i}} 		& \le & \norm{\sum_{i=1}^k a_i x_{n_i}}+\norm{\sum_{i=1}^\ell a_i x_{n_i}}\\	
		& \le & A_s\norm{\sum_{i=1}^k a_i x_{i}}+ A_s\norm{\sum_{i=1}^\ell a_i x_{i}} \\
		& \le & A_s(1+C_b)\norm{\sum_{i=1}^k a_i x_i}=A_s(1+C_b)N(\1),
\end{eqnarray*}

it follows that 

\begin{equation*}
\max \{\norm{\sum_{i=1}^\ell a_i x_{n_{\pi(i)}}}, \norm{\sum_{i=\ell+1}^k a_i x_{n_{\pi(i)}}}, \norm{\sum_{i=\ell+1}^k a_i x_{n_i}} , \norm{\sum_{i=1}^\ell a_i x_{n_i}}\}\le A_s(1+C_b)N(\1).
\end{equation*}
If we define for $(m_1,\dots,m_\ell)\in\Nk$ 
$$f(m_1,\dots,m_\ell)=\frac{1}{A_s(1+C_b)N(\1)} \sum_{i=1}^\ell a_i x_{m_i}$$

and for $(m_1,\dots,m_{k-\ell})\in\Sk{\bN}{k-\ell}$ 

$$g(m_{1},\dots,m_{k-\ell})=\frac{1}{A_s(1+C_b)N(\1)} \sum_{i=1}^{k-\ell} a_{\ell+i} x_{m_i},$$

then $f$ and $g$ take values in the unit ball of $\banX$. Observe now that

\begin{eqnarray*}
N(\1)\le B_s\norm{\sum_{i=1}^k a_i x_{n_i}} & \le & B_s\norm{-\sum_{i=1}^\ell a_i x_{n_i}+\sum_{i=\ell+1}^k a_i x_{n_i}}+2B_s\norm{\sum_{i=1}^\ell a_i x_{n_i}}\\
		& \le & B_s\norm{\sum_{i=1}^\ell a_i x_{n_i}-\sum_{i=\ell+1}^k a_i x_{n_i}}+2B_sC_b\norm{\sum_{i=1}^\ell a_i x_{n_i}-\sum_{i=\ell+1}^k a_i x_{n_i}}\\
		& = & B_s(2C_b+1)\norm{\sum_{i=1}^\ell a_i x_{n_i}-\sum_{i=\ell+1}^k a_i x_{n_i}},
\end{eqnarray*}

and hence

\begin{equation*}
\inf_{\nbar\in\Nk} \dX(f(n_1,\dots,n_\ell),g(n_{\ell+1},\dots,n_k)) \ge \frac{1}{B_s(2C_b+1)A_s(1+C_b)}.
\end{equation*}

If $C$ is any number satisfying $\omega(C)\le \frac{1}{2K_u}\rho( \frac{1}{B_s(2C_b+1)A_s(1+C_b)})$  it follows from Lemma \ref{lem:metric-uniform} that

\begin{equation*}
\sup_{\nbar\in\Nk} \norm{f(n_{\pi(1)},\dots,n_{\pi(\ell)})-g(n_{\pi(\ell+1)},\dots,n_{\pi(k)})}\ge C,
\end{equation*}

and thus there exists $\nbar\in \Nk$ such that 

\begin{equation*}
N(\vep)=\norm{\sum_{i=1}^\ell a_i x_{n_{\pi(i)}}-\sum_{i=\ell+1}^k a_i x_{n_{\pi(i)}}}\ge \frac{C}{2}A_s(1+C_b)N(\1).
\end{equation*}

Taking $\gamma=\max\{1,\frac{C}{2}A_s(1+C_b)\}$ in \eqref{eq:uncond} concludes the proof.
\end{proof}

The next theorem shows that upper stability is a strengthening of Kalton's property $\cQ$. 

\begin{theo}\label{thm:stable-Q}
Every upper stable metric space has property $\cQ$.
\end{theo}

\begin{proof}
Let $\metXd$ be $K_u$-upper stable and fix $k\ge 1$ and a Lipschitz map $f\colon (\Nk,\sd_\sI)\to \metX$.
Define 
$$\cA\eqd\{\nbar\in \Sk{\bN}{2k}\colon \dX(f(n_1,\dots,n_k),f(n_{k+1},\dots, n_{2k})\le 2K_u\cdot\lip(f)\}.$$
By Ramsey theorem there exists an infinite subset $\bM$ of $\bN$ such that either $[\bM]^{2k}\subset \cA$ or $[\bM]^{2k}\cap \cA=\emptyset$.
If the first possibility happens then for all $\nbar\in [\bM]^{2k}$ we have 
$$\dX(f(n_1,\dots,n_k),f(n_{k+1},\dots, n_{2k}))\le 2K_u\cdot\lip(f).$$
If $\nbar,\mbar\in[\bM]^k$ we can choose $\xbar\in[\bM]^k$ such that $x_1>\max\{n_k,m_k\}$ and thus
$$\dX(f(\nbar),f(\mbar))\le \dX(f(\nbar),f(\xbar))+\dX(f(\xbar),f(\mbar))\le 4K_u\cdot\lip(f).$$
It remains to show that the second possibility cannot happen. Consider the interlacing permutation $\pi\colon [2k]\to [2k]$ defined by 
\begin{equation*}
\pi(i)=\begin{cases}
2i-1 \quad if \quad 1\le i\le k\\
2(i-k) \quad if \quad k+1\le i\le 2k.
\end{cases}
\end{equation*}

It is immediate to verify that $\pi$ preserves the order on $\{1,\dots,k\}$ and on $\{k+1,\dots,2k\}$ and if $\nbar\in[\bM]^{2k}$ then

\begin{equation*}
\dX(f(n_{\pi(1)},\dots,n_{\pi(k)}),f(n_{\pi(k+1)},\dots,n_{\pi(2k)}))=\dX(f(n_1,n_3,\dots,n_{2k-1}),f(n_2,n_4,\dots,n_{2k}))\le \lip(f).
\end{equation*}
Therefore,

\begin{equation*}
\sup_{\nbar\in [\bM]^{2k}}\dX(f(n_{\pi(1)},\dots,n_{\pi(k)}),f(n_{\pi(k+1)},\dots,n_{\pi(2k)}))\le \lip(f),
\end{equation*}

and by $K_u$-upper stability it holds

\begin{equation*}
\inf_{\nbar\in [\bM]^{2k}}\dX(f(n_{1},\dots,n_{k}),f(n_{k+1},\dots,n_{2k}))\le K_u\lip(f).
\end{equation*}

In particular there exists $\nbar\in[\bM]^{2k}$ such that 

\begin{equation*}
\dX(f(n_{1},\dots,n_{k}),f(n_{k+1},\dots,n_{2k}))\le 2K_u\lip(f),
\end{equation*}

but this implies that $[\bM]^{2k}\cap \cA\neq\emptyset$; a contradiction.
\end{proof}

Recall that a Banach space $\banX$ has trivial Rademacher type if and only if $\banX$ contains the $\ell_1^n$'s (c.f. \cite{MaureyPisier76}). It was shown by Guerre-Delabriere and Laprest\'e \cite{GL81} that every non-reflexive stable Banach space contains an isomorphic copy of $\ell_1$ and hence has trivial Rademacher type. Raynaud showed that every non-reflexive Banach space whose unit ball uniformly embeds into a stable metric has a spreading model isomorphic to $\ell_1$. Without delving too deep into the theory of spreading models, a spreading model of $\banX$ is a Banach space that can be associated to $\banX$ and that is finitely representable in $\banX$. Therefore the unit ball of a non-reflexive Banach space with non-trivial Rademacher type does not even embed uniformly into a stable space. Classical examples of non-reflexive Banach space with non-trivial Rademacher type are James non-reflexive space of type $2$ \cite{James78} and Pisier-Xu interpolation spaces \cite{PisierXu87}. 
Raynaud's result was extended by Kalton. Indeed, it is a result from \cite{Kalton07} that if $\banX$ has property $\cQ$ and is non-reflexive then $\banX$ has a spreading model isomorphic to $\ell_1$. Kalton's proof relies on a result of Beauzamy \cite{Beauzamy79} which says that a Banach space has the alternating Banach-Saks property if and only if none of its spreading models are isomorphic to $\ell_1$. Recall that a Banach space $\banX$ has the alternating Banach-Saks property if every bounded sequence $\xn\in \banX$ has a subsequence $\yn$ so that the alternating Cesaro means $\frac{1}{n}\sum_{k=1}^{n}(-1)^k y_k$ converge in norm to $0$. Note that the class of Banach spaces with the alternating Banach-Saks property strictly extends the class of Banach spaces with non-trivial type since $\co$ has the alternating Banach-Saks property and trivial Rademacher type.

The conclusion of Raynaud and Kalton results obviously holds for Since upper stability implies property $\cQ$ it follows from Kalton's result that a non-reflexive and upper stable Banach space has a spreading model isomorphic to $\ell_1$, and thus contains the 
$\ell_1^n$'s. If we are merely interested in the containment of the $\ell_1^n$'s  the elementary argument given below is all we need.

\begin{theo}\label{thm:type-stable-reflexive}
Let $\banX$ be a Banach space which is non-reflexive. Assume that, either
\begin{enumerate}[(i)]
\item $\banX$ coarsely, or uniformly, embeds into an upper stable metric space
\item[] or
\item the unit ball of $\banX$ uniformly embeds into an upper stable metric space.
\end{enumerate}
Then $\banX$ contains the $\ell_1^n$'s.
\end{theo}

\begin{proof}
Assertion $(i)$ is a simple consequence of Lemma \ref{lem:coarse-upper}. Since $\banX$ is not reflexive there exists a normalized James' sequence $\xn$ such that for all $1\le \ell\le k$, $b_1,\dots,b_k\ge 0$, and $\nbar\in\Nk$,
\begin{equation*}
\norm{\sum_{i=1}^\ell b_i x_{n_i}-\sum_{i=\ell+1}^{k}b_i x_{n_i}}\ge \sum_{i=1}^{k}b_i
\end{equation*}

Let $a_1,\dots,a_k\in\bR$ and choose $\vep_1,\dots,\vep_k\in\{\pm1\}$ such that $\vep_i\abs{a_i}=a_i$ for all $1\le i\le k$. Let $\ell\eqd|\{i\colon \vep_i=+1\}|$ and $\pi$ the permutation on $[k]$ that maps $\{1,\dots,\ell\}$ onto $\{i\colon \vep_i=+1\}$ and $\{\ell+1,\dots,k\}$ onto $\{i\colon \vep_i=-1\}$ while preserving the order on the respective sets. Let $\alpha\eqd \sum_{i=1}^{k}\abs{a_i}$ and define for $\mbar\in\Nk$, $f(m_1,\dots,m_\ell)=\frac{1}{\alpha}\sum_{i=1}^\ell \abs{a_i} x_{m_i}$ and for $\mbar\in\Sk{\bN}{k-\ell}$, $g(m_{1},\dots,m_{k-\ell})= \frac{1}{\alpha}\sum_{i=1}^{k-\ell} \abs{a_{\ell+i}} x_{m_i}$. It follows from Lemma \ref{lem:coarse-upper} that
\begin{equation*}
\inf_{\nbar\in\Nk} \norm{f(n_1,\dots,n_\ell)-g(n_{\ell+1},\dots,n_k)}\le C \sup_{\nbar\in\Nk} \norm{f(n_{\pi(1)},\dots,n_{\pi(\ell)})-g(n_{\pi(\ell+1)},\dots,n_{\pi(k)})}.
\end{equation*}

Therefore, for all $\nbar\in\Nk$

\begin{equation*}
\frac{1}{C}\sum_{i=1}^{k}\abs{a_i}\le \sup_{\nbar\in\Nk}\norm{\sum_{i=1}^k a_i x_{n_i}},
\end{equation*}
and there exists $\nbar\in\Nk$ such that
\begin{equation*}
\frac{1}{2C}\sum_{i=1}^{k}\abs{a_i}\le \norm{\sum_{i=1}^k a_i x_{n_i}}\le \sum_{i=1}^{k}\abs{a_i}.
\end{equation*} 

For assertion $(ii)$ we can apply Lemma \ref{lem:metric-uniform} since $f$ and $g$ take values in the unit ball of $\banX$.
\end{proof}
\begin{coro}\label{cor:type-stable-reflexive}
Let $\banX$ be a Banach space that does not contain the $\ell_1^n$'s. Assume that, either
\begin{enumerate}[(i)]
\item $\banX$ coarsely, or uniformly, embeds into an upper stable metric space
\item[] or
\item the unit ball of $\banX$ uniformly embeds into a upper stable metric space.
\end{enumerate}
Then $\banX$ is is reflexive.
\end{coro}

The assumptions in Theorem \ref{thm:type-stable-reflexive} and Corollary \ref{cor:type-stable-reflexive} cannot be lifted since $\ell_1$ is stable, contains the $\ell_1^n$'s and is not reflexive. 

\section{Final remarks and open problems}

\subsection{Optimal factor in Ostrovskii's determinacy theorem}
Ostrovskii's finite determinacy theorem says that there is a universal constant $\alpha\in(0,\infty)$ such that for every Banach space $\banY$ and every locally finite metric space $\metX$, $\metX$ admits a bi-Lipschitz embedding into $\banY$ with distortion at most $\alpha\cdot \beta$ whenever $\metX$ is finitely $\beta$-Lipschitzly representable in $\banY$.  A quick inspection of Ostrovskii's original proof gives that we can take $\alpha=2100$. The study of the parameter $\alpha$ was studied in \cite{OO19a} and \cite{OO19b}. Most notably it was shown that in general necessarily $\alpha>1$ and if the host space is an $\ell_p$ -sum of nested finite-dimensional space then $\alpha$ can be taken to be arbitrarily close to ~1. Logarithmic spiral gluing was introduced to achieve this degree of precision in this special case.

\subsection{$L_p[0,1]$-compression exponent vs $\ell_p$-compression exponent}\label{sec:compression}
The \emph{compression exponent} of $\metXd$ in $\metYd$ ($\metY$-compression of $\metX$ in short) introduced by Guentner and Kaminker \cite{GuentnerKaminker04}, is the parameter denoted by $\alpha_\metY(\metX)$ and defined as the supremum of all numbers $0\le \alpha\le 1$ for which there exist $f\colon \metX\to\metY$, $\tau\in(0,\infty)$, and $A\in[1,\infty)$ such that $\dX(x,y)\ge \tau$ implies
$$\frac{1}{A}\dX(x,y)^{\alpha}\le \dY(f(x),f(y))\le A \dX(x,y).$$
The compression exponent of an (unbounded) metric space is clearly invariant under coarse bi-Lipschitz embeddings. We will write $L_p$ for $L_p[0,1]$. For any metric space $\metX$, since $\ell_p$ embeds isometrically into $L_p$, and since $L_p$ contains an isometric copy of $L_2$, the inequalities $\alpha_{\ell_p}(\metX)\le \alpha_{L_p}(\metX)$ and $\alpha_2(\metX)\le \alpha_{L_p}(\metX)$ always hold for any metric space $\metX$. It is well known that finite subsets of $L_p$ embed isometrically into $\ell_p$ and thus it follows that the $L_p$-compression of an unbounded sequence of finite metric spaces (defined in a natural way) and its $\ell_p$-compression coincides. A similar reasoning does not work to prove an analogue statement if one is dealing with an (unbounded) infinite metric space since, for instance, $L_4$ does not coarse bi-Lipschitz embeds into $\ell_4$ \cite{KaltonLova08}. Corollary \ref{cor:locfincomplp}, which obviously applies to finitely generated groups, answers negatively a question raised by Naor and . According to Naor and Peres this subtlety between the $L_p$-compression and the $\ell_p$-compression for infinite groups was first pointed out by Marc Bourdon, and Naor and Peres  asked (Question 10.7 in \cite{NaorPeres11}) whether there is a finitely generated group $\Gamma$ such that $\alpha_{\ell_p}(\Gamma)\neq \alpha_{L_p}(\Gamma)$. Theorem \ref{thm:Lp} below is an easy consequence of Theorem \ref{thm:barycentric-internal}.
  
\begin{theo}\label{thm:Lp}
Let $p\in[1,\infty]$. If $\ell_p$ is finitely representable in a Banach space $\banY$ then any locally finite subset of $L_p$ admits a bi-Lipschitz embedding into $\banY$. In particular,  any locally finite subset of $L_p$ admits a bi-Lipschitz embedding into $\ell_p$. 
\end{theo}

In particular, Theorem \ref{thm:Lp} and the coarse Lipschitz invariance of compression exponents provides a complete answer for Naor and Peres question.
 
\begin{coro}\label{cor:locfincomplp}
Let $\metX$ be a proper metric space. Then for every $p\in[1,+\infty]$
\begin{equation*}
\alpha_{\ell_p}(\metX)=\alpha_{L_p}(\metX).
\end{equation*}
In particular, if $\Gamma$ is a finitely generated group with its canonical word metric, or a compactly generated group with its canonical proper metric, then 
\begin{equation*}
\alpha_{\ell_p}(\Gamma)=\alpha_{L_p}(\Gamma).
\end{equation*}
\end{coro}

\begin{rema}
Corollary \ref{cor:locfincomplp} also applies to locally compact second countable groups by the results of \cite{HP06}.
\end{rema}

The cases $p=\infty$ and $p=2$ of Theorem \ref{thm:Lp} are just reformulations of Corollary \ref{cor:BL} and Corollary \ref{cor:Hilbert}, respectively. When $1<p<\infty$ simply observe that every finite subset of $L_p$ embeds isometrically into some $\ell_p^n$, and thus $L_p$ is finitely Lipschitzly representable into every Banach space which contains the $\ell_p^n$'s.

A coarse version of Theorem \ref{thm:Lp} seems to have appeared first in the unpublished manuscript \cite{Ostrovskii_unpublished}. After the appearance of \cite{BaudierLancien08}, the coarse statement of \cite{Ostrovskii_unpublished} was upgraded to the bi-Lipschitz category, and only the explicit statement in Corollary \ref{cor:Hilbert} made it to the published paper \cite[Theorem 6]{Ostrovskii09} while the explicit statement of Theorem \ref{thm:Lp} was surprisingly left out (but can be derived from \cite[Proposition 1]{Ostrovskii09}). Theorem \ref{thm:Lp} appeared explicitely in \cite{Baudier12} where the statement and the proof uses the convenient framework of $\mathcal{L}_{p}$-spaces. 

Naor and Peres \cite{NaorPeres08} proved that $\alpha^{\#}_2(\Gamma)\le \alpha^{\#}_{L_p}(\Gamma)$, for any $p\in[1,\infty)$ and every finitely generated amenable group $\Gamma$, where $\alpha^{\#}_\metY$ denotes the $\metY$-equivariant compression (see \cite{NaorPeres08} for the definition). The following result follows from Theorem \ref{thm:cfrLipschitz}.

\begin{coro}\label{cor:proper-comp}
Let $\metX$ be a proper metric space. 
If $\banX$ and $\banY$ are two Banach spaces such that $\banX$ is finitely representable in $\banY$, then $\alpha_\banX(\metX)\le \alpha_\banY(\metX)$. In particular, $\alpha_2(\metX)\le \alpha_\banY(\metX)$ for every infinite-dimensional Banach space $\banY$.
\end{coro} 
Note that the second part of the statement can be obtained with Corollary \ref{cor:Hilbert} only. We do not know if an equivariant analogue of Corollary \ref{cor:proper-comp} is true.
\begin{prob}
Let $\Gamma$ be a finitely generated group. Do we have $\alpha^{\#}_2(\Gamma)\le \alpha^{\#}_\banY(\Gamma)$ for every infinite-dimensional Banach space $\banY$?
\end{prob} 

The equivariant analogue of Corollary \ref{cor:locfincomplp} also seems open.

\begin{prob}
Let $\Gamma$ be a finitely generated group. Do we have for every $p\in[1,+\infty)$ $$\alpha^{\#}_{\ell_p}(\Gamma)=\alpha^{\#}_{L_p}(\Gamma)?$$
\end{prob} 

\subsection{Equivariant versions of the embedding results}
Brown and Guentner showed that every countable discrete group admits a proper affine isometric action on the reflexive space $\xbg$. This result was extended to locally compact second countable groups in \cite{HP06}. A coarse equivariant version of Kalton's embedding result for stable metrics was proved by Rosendal \cite[Theorem 22]{Rosendal17}. It would be interesting to know if there are equivariant versions of Corollary \ref{cor:BL} for countable discrete groups and Corollary \ref{cor:proper-CL} and Corollary \ref{cor:proper-cotype} for compactly generated groups.

\subsection{Embeddability of reflexive spaces into stable space}

Note that $T$ does not have the alternating Banach-Saks property, in a strong sense, since all its spreading models are isomorphic to $\ell_1$. However there are reflexive Banach spaces with the alternating Banach-Saks property which fail to be stable. Problem \ref{pb:Kalton} which was raised by Kalton \cite[Problem 6.1-6.2]{Kalton07} asks for a converse to Theorem \ref{thm:Kalton}.

\begin{prob}\label{pb:Kalton}
Let $X$ be a separable reflexive Banach space. 
\begin{enumerate}[(i)]
\item Does $X$ coarsely embeds into a stable metric space?
\item Does $X$ uniformly embeds into a stable metric space?
\item Does the unit ball of $X$ uniformly embeds into a stable metric space?
\end{enumerate}
\end{prob}


The difficulty of Problem \ref{pb:Kalton} can be partially explained by the deep work from \cite{Kalton07} about the relationship between reflexivity and property $\cQ$. Firstly, one of the main results from \cite{Kalton07} is that every reflexive Banach space has property $\cQ$ and this was used to resolve negatively the long-standing problem whether $\co$ (the unit ball of $\co$) coarsely (uniformly) embeds into a reflexive Banach space. Secondly, if $\banX$ is a Banach space with the alternating Banach-Saks property then under either assumption of Problem \ref{pb:Kalton}, $\banX$ must necessarily be reflexive. We suspect the answer to assertion $(i)$ and $(ii)$ in Problem \ref{pb:Kalton} are negative since we conjecture that there are reflexive spaces that are not upper stable. In particular, we conjecture that upper stability is strictly stronger than property $\cQ$.

Note that it follows from \cite{Raynaud83} that the unit ball of a Tsirelson-like space does not uniformly embed into a super-stable space. It was shown in \cite{BragaSwift19} that  $\ban T^*$ (Tsirelson's original space that contains the $\ell_\infty^n$'s which is the dual of Figiel-Johnson-Tsirelson space $\ban T$ \cite{FigielJohnson74}) does not coarsely embeds into a super-stable space. However it is not known if this is true for all Tsirelson-like spaces since the following problem is still open (\cite[Problem 6.6]{Kalton07}).

\begin{prob}
If $\banX$ coarsely embeds into a super-stable space, does $\banX$ contain an isomorphic copy of $\ell_p$ for some $p\in[1,\infty)$?
\end{prob}

\subsection{Reflexive asymptotic-$\ell_2$ spaces, Baum-Connes and Novikov conjecture}
 
To the best of our knowledge there is no example yet of bounded geometry spaces failing the coarse Baum-Connes conjecture or the coarse Novikov conjecture. In \cite{Yu00} and \cite{KasparovYu06} the conjectures are proved for bounded geometry spaces admitting a weak embedding (coarse embedding) into very regular Banach spaces (Hilbert space or a uniformly convex Banach space with a uniformly convex dual). It would be interesting to investigate whether the conjectures can be proven under a much stronger embeddability requirement (bi-Lipschitz embedding) into Banach spaces with weaker geometric features. For instance, the space $\xbl$ is reflexive, asymptotic-$\ell_2$, asymptotically uniformly convex with an asymptotically uniformly convex dual. 
 
\begin{prob}\label{pb:novikov}
Let $\metX$ be a discrete metric space with bounded geometry. Assume that $\metX$ admits a bi-Lipschitz embedding into a reflexive asymptotic-$\ell_2$ Banach space. Does $\metX$ satisfy the coarse Baum-Connes or the coarse geometric Novikov conjecture?
\end{prob}

An important warning is in order here. By Corollary \ref{cor:BL} a positive answer to Problem \ref{pb:novikov} will imply a positive answer to the respective conjecture for every metric space with bounded geometry! Nevertheless, it seems important to have a deeper understanding of the tradeoff between the faithfulness of the embedding and the geometric features of the host space in order to be able to confirm the conjectures. 

\subsection{More on stability}
Kalton \cite{Kalton07} proved that James space $\ban J$ and its dual $\ban J^*$ \cite{James50,James51} do not have the $\cQ$-property, and this gives two more examples of non-reflexive spaces that do not embed coarsely into a stable space.

The work of Rosendal \cite{Rosendal17} (and the references therein) contains lots of information about the existence of compatible stable metrics on groups.

To find more examples of Banach spaces which are not stable we need to resort to non-commutative $L_p$-spaces. We refer to \cite{PisierXu_HB03} for a thorough discussion of non-commutative $L_p$-spaces. The commutative $L_p$-spaces belong to the class of non-commutative $L_p$-spaces. For instance, $L_p[0,1]$ can been represented as the non-commutative $L_p$-space associated to $L_\infty[0,1]$ considered as a von Neuman algebra on the Hilbert space $L_2[0,1]$. The simplest truly non-commutative $L_p$-spaces are the Schatten classes (denoted $\ban{S}_p(H)$). They are defined as the non-commutative $L_p$-spaces associated to $B(H)$, the algebra of all bounded operators on a Hilbert space $H$, equipped with be the usual trace on $B(H)$. The commutative theory can be satisfactorily extended to a large extent to the non-commutative setting, however there are some significant differences. The stability property of the non-commutative spaces is one of them.
\begin{theo}
Let $p\in[1,\infty)$, $p\neq 2$. Then a non-commutative $L_p$-space associated to a von Neumann algebra $\mathfrak{M}$ is stable if and only if $\mathfrak{M}$ is of type \textrm{I}.
\end{theo} 
The ``if'' part of the theorem above was independently proved in \cite{Arazy83}, and \cite{Raynaud84a} where it is shown that the non-commutative $L_p$-space associated to a von Neumann algebra $\mathfrak{M}$ of type $\textrm{I}$ can be written as an $\ell_p$-sum of commutative vector-valued $L_p$-spaces whose values fall into stable Banach spaces. The ``only if'' part comes from \cite{Marcolino97} and is proved in two steps as follows. Marcolino first showed that if $\mathfrak{M}$ is a von Neuman algebra not of type $\textrm{I}$, then for $1\le p\le \infty$, the non-commutative $L_p$-space associated to the hyper finite $\textrm{II}_1$ factor is isometric to a (1-complemented) subspace of the non-commutative $L_p$-space associated to $\mathfrak{M}$. The conclusion follows from the fact that the non-commutative $L_p$-space associated to the hyper finite $\textrm{II}_1$ factor is not stable. Those spaces have a completely different linear structure compared to Tsirelson-like spaces since they contains copies of $\ell_p$. The following problem seems open.

\begin{prob}
Does the non-commutative $L_p$-space associated to the hyper finite $\textrm{II}_1$ factor admit an equivalent stable norm?
\end{prob}

\bibliographystyle{alpha}



\end{document}